\documentclass[reqno]{amsart}
\pdfoutput=1
\usepackage{amssymb}
\usepackage{graphicx}
\usepackage{mathrsfs}
\usepackage{hyperref}
\usepackage{amsrefs}

\theoremstyle{plain}
    \newtheorem{theorem}{Theorem}[section]

\theoremstyle{definition}
    \newtheorem{definition}[theorem]{Definition}

\theoremstyle{remark}
    \newtheorem{remark}[theorem]{Remark}

    \newcommand{\IGNORE}[1]{}

    \newcommand{\B}{{\mathrm{b}}}
    \newcommand{\R}{{\mathbb{R}}}           
    \newcommand{\U}{{\mathbb{U}}}           

    \newcommand{\mb}[1]{\mathbf{#1}}

    \newcommand{\BP}{{\mathbf{p}}}          
    \newcommand{\BR}{{\mathbf{r}}}          
    \newcommand{\BT}{{\mathbf{t}}}          
    \newcommand{\BU}{{\mathbf{u}}}          
    \newcommand{\BX}{{\mathbf{x}}}          
    \newcommand{\BZ}{{\mathbf{z}}}          

    \newcommand{\GS}{\geqslant}
    \newcommand{\HA}{{\textstyle\frac{1}{2}}}
    \newcommand{\ID}{{\mathrm{id}}}
    
    \newcommand{\LS}{\leqslant}

    \newcommand{\C}{{\mathscr{C}}}          
    \newcommand{\D}{{\mathscr{D}}}          
    \renewcommand{\L}{{\mathscr{L}}}        
    \newcommand{\CB}{{\C_\B}}               
    \newcommand{\SP}{{\mathscr{P}}}         

    \newcommand{\DST}{\displaystyle}
    
    \newcommand{\EPS}{\varepsilon}
    
    \newcommand{\IND}{{\mathbf{1}}}
    \newcommand{\INT}{{\mathcal{U}}}
    \newcommand{\mc}[1]{\mathcal{#1}}
    \newcommand{\LEB}{{\mathfrak{L}}}

    \newcommand{\RHO}{\varrho}
    \newcommand{\TAN}{{\mathbb{T}}}
    \newcommand{\TST}{\textstyle}
    
    \newcommand{\WAS}{{\mathbf{W}}}

    \newcommand{\jt}{\textstyle}
    \newcommand{\jd}{\displaystyle}
    \newcommand{\e}[1]{{(#1)}}
    \newcommand{\opn}[1]{\operatorname{#1}}

    \DeclareMathOperator*{\ARGMIN}{{\mathrm{argmin}}}

    \numberwithin{equation}{section}

\title
    [Variational Particle Scheme]
    {Variational Particle Schemes for the Porous Medium Equation
     and for the System of Isentropic Euler Equations}

\author
    {Michael Westdickenberg and Jon Wilkening}

\address
    {Michael Westdickenberg,
     School of Mathematics,
     Georgia Institute of Technology,
     686 Cherry Street,
     Atlanta, GA 30332-0160,
     U.S.A.}
\email
    {mwest@math.gatech.edu}

\address
    {Jon Wilkening,
     Department of Mathematics,
     University of California,
     1091 Evans Hall \#3840,
     Berkeley, CA 94720-3840,
     U.S.A.}
\email
    {wilken@math.berkeley.edu}

\date{\today}

\subjclass[2000]
   {35L65, 49J40, 76M30, 76M28}

\keywords
   {Optimal Transport, Wasserstein Metric,
    Isentropic Euler Equations,
    Porous Medium Equation, Numerical Methods}

\begin{document}

\begin{abstract}
  Both the porous medium equation and the system of isentropic Euler
  equations can be considered as steepest descents on suitable
  manifolds of probability measures in the framework of optimal
  transport theory. By discretizing these variational
  characterizations instead of the partial differential equations
  themselves, we obtain new schemes with remarkable stability
  properties. We show that they capture successfully the nonlinear
  features of the flows, such as shocks and rarefaction waves for the
  isentropic Euler equations. We also show how to design higher order
  methods for these problems in the optimal transport setting using
  backward differentiation formula (BDF) multi-step methods or
  diagonally implicit Runge-Kutta methods.
\end{abstract}

    \maketitle
    \tableofcontents

\section{Introduction}\label{S:I}

It is well-known that the dynamics of many physical systems can be
derived from variational principles. For classical mechanics, for
example, the first variation of an {\em action functional} yields the
equations of motion, which are also called the Euler-Lagrange
equations of the functional. By discretizing these equations one can
obtain numerical methods. This strategy is time-tested and reliable.
More recently, however, an alternative approach has
attracted considerable interest, which consists of {\em discretizing
  the action functional} instead. The numerical method is then given
by the Euler-Lagrange equations of this discrete functional. Methods
obtained in this way enjoy remarkable stability properties and
preserve the ``symplectic'' structure of the problem. We refer the
reader to \cite{MarsdenWest} for further information. In the present
paper, we apply a similar strategy for the porous medium equation and
the system of isentropic Euler equations in one space dimension, and
derive numerical methods by discretizing a suitable variational
principle.

The isentropic Euler equations model the dynamics of compressible fluids
under the simplifying assumption that the thermodynamical entropy is
constant in space and time. These equations form a system of
hyperbolic conservation laws for the density $\RHO$ and the Eulerian
velocity field $\BU$. They take the form
\begin{equation}
\begin{aligned}
    & \partial_t \RHO + \nabla\cdot(\RHO\BU) = 0,
\\
    & \partial_t ( \RHO\BU ) + \nabla\cdot( \RHO\BU\otimes\BU )
        + \nabla P(\RHO) = 0,
\end{aligned}
\label{E:IEE}
\end{equation}
where $(\RHO,\BU) \colon [0,\infty)\times\R^d \longrightarrow \U$ are
measurable functions with $\U := [0,\infty)\times\R^d$. Since entropy
is assumed to be constant, the pressure $P(\RHO)$ depends on the
density only. It is defined in terms of the internal energy $U(\RHO)$
of the fluid by
\begin{equation}
    P(\RHO) := U'(\RHO)\RHO - U(\RHO)
    \quad\text{for all $\RHO\GS 0$.}
\label{E:PRESS}
\end{equation}
%
For the important case of polytropic fluids, we have
\begin{equation}
    U(\RHO) = \frac{\kappa\RHO^\gamma}{\gamma-1}
    \quad\text{and}\quad
    P(\RHO)=\kappa\RHO^\gamma.
\label{E:POLYT}
\end{equation}
Here $\gamma>1$ is the adiabatic coefficient, and we will assume the
common normalization $\kappa := \theta^2/\gamma$ with $\theta :=
(\gamma-1)/2$.
For isothermal flows, we have
\begin{equation}
    U(\RHO) = \RHO\log\RHO
    \quad\text{and}\quad
    P(\RHO) = \RHO.
\label{E:ISO}
\end{equation}
We are interested in the Cauchy problem, so we require that
$$
    (\RHO,\BU)(t=0,\cdot) = (\bar{\RHO},\bar{\BU})
$$
for suitable initial data $(\bar{\RHO},\bar{\BU})$ with finite mass and total energy.

It is again well-known \cites{ArnoldKhesin, HolmMarsdenRatiu} that the
system of isentropic Euler equations can be considered as the
Euler-Lagrange equations of a certain action functional, so we might try to derive numerical schemes by following the procedure outlined
above. But we want to work in the framework of optimal transport theory.
Our discussion in the present paper will therefore not be based on the variation of an action functional, but on the notion that the system of isentropic Euler equations \eqref{E:IEE} can be considered as a {\em steepest descent} on a suitable abstract manifold. This notion was introduced in \cite{GangboWestdickenberg} and will be explained in detail below. To motivate our approach, recall that any {\em smooth} solution $(\RHO,\BU)$ of \eqref{E:IEE} also satisfies
\begin{equation}
    \partial_t \bigg( \frac{1}{2}\RHO|\BU|^2+U(\RHO) \bigg)
        + \nabla\cdot \Bigg( \bigg( \frac{1}{2}\RHO|\BU|^2
        + Q(\RHO) \bigg) \BU \Bigg) = 0,
\label{E:IEE2}
\end{equation}
where $Q(\RHO) := U(\RHO) + P(\RHO)$. Therefore the functions
$$
    \eta(\RHO,\BU) := \frac{1}{2}\RHO|\BU|^2 + U(\RHO)
    \quad\text{and}\quad
    q(\RHO,\BU) := \bigg( \frac{1}{2}\RHO|\BU|^2 +
        Q(\RHO) \bigg) \BU
$$
form an entropy/entropy flux pair for \eqref{E:IEE}; see
\cite{Dafermos}. The entropy $\eta(\RHO, \BU)$ is a convex function of
$\RHO$ and $\RHO\BU$. As a consequence of \eqref{E:IEE2}, the total
energy
\begin{equation}
    E[\RHO,\BU](t) := \int_{\R^d}
        \Big( \HA\RHO|\BU|^2 + U(\RHO) \Big)(t,x) \,dx
\label{E:TOTL}
\end{equation}
is conserved in time as long as the solution remains smooth.

On the other hand, it is known that solutions of \eqref{E:IEE} are
generally not smooth: No matter how regular the initial data is, jump
discontinuities can form in finite time. These jumps occur along
$d$-dimensional submanifolds in space-time and are called shocks.
Across shocks, total energy is dissipated (that is, transformed into
forms of energy such as heat that are not accounted for in an
isentropic model), so \eqref{E:IEE2} cannot hold anymore. It is
therefore natural to consider {\em weak} solutions $(\RHO,\BU)$ of
\eqref{E:IEE} such that $E[\RHO,\BU]$ is a non-increasing function of
time.

As will be explained below, we consider solutions of the isentropic
Euler equations as curves on the space of probability measures, using optimal transport theory. We are interested in curves along which the total energy decreases as fast as possible, under conservation of mass and momentum. This steepest descent interpretation has
already proved very fruitful for certain degenerate parabolic
equations, such as the porous medium equation; see below. Numerical
schemes based on this variational principle have been derived in
\cite{KinderlehrerWalkington}. As the authors point out, these schemes
have remarkable properties: oscillations can be reduced because the optimizations take place in weak topologies; the approximations can be constructed using discontinuous functions because no derivatives are
needed; and overall, the schemes are very stable.
In this paper, we revisit the discretization of the porous medium
equation, and also introduce variational approximations for the
isentropic Euler equations.  We demonstrate that our schemes capture
very well the nonlinear features of the flows.  The simplest versions
of our schemes are first order accurate.

Designing higher order methods in the optimal transport setting is an
interesting challenge that has not previously been addressed.
We present schemes that achieve second order
convergence (away from shocks) and propose a general strategy for
constructing higher order methods in a backward
differentiation formula multi-step (BDF) or diagonally
implicit Runge-Kutta framework. Again we observe that our schemes,
due to their variational nature, are remarkably stable.
\medskip

To put our work and that in \cite{GangboWestdickenberg} in
perspective, let us first recall recent research by various authors on
a steepest descent interpretation for certain degenerate parabolic
equations. It was shown by Otto \cite{Otto} that the porous medium
equation
\begin{equation}
    \partial_t\RHO - \Delta P(\RHO) = 0
    \quad\text{in $[0,\infty)\times\R^d$}
\label{E:POR}
\end{equation}
is a {\em gradient flow} in the sense enumerated below. This result
has been generalized considerably, and we refer the reader to
\cites{AmbrosioGigliSavare, Villani} for further information.

\medskip

(1) We denote by $\SP(\R^d)$ the space of all $\LEB^d$-measurable,
non-negative functions with unit integral and finite second
moments, where $\LEB^d$ is the Lebesgue measure. The space
$\SP(\R^d)$ is equipped with the {\em Wasserstein distance},
defined by
\begin{equation}
    \WAS(\RHO_1,\RHO_2)^2
        := \inf\Bigg\{ \iint_{\R^d\times\R^d}
                |x_2-x_1|^2 \,\gamma(dx_1,dx_2) \colon
            \pi^i\#\gamma = \RHO_i\LEB^d \Bigg\},
\label{E:WAS}
\end{equation}
where the infimum is taken over all transport plans $\gamma$, which are probability measures on the product space $\R^d\times\R^d$ with the property that the pushforward $\pi^i\#\gamma$ of $\gamma$ under the projection $\pi^i:\R^d\times\R^d\rightarrow\R^d$ onto the $i$th component equals $\RHO_i\LEB^d$. The infimum in \eqref{E:WAS} represents the minimal quadratic cost required to transport the measure $\RHO_1\LEB^d$ to the measure
$\RHO_2\LEB^d$, and one can show that it is always attained for some optimal transport plan $\gamma$. In fact, there is a lower semicontinuous, convex function $\phi\colon\R^d\longrightarrow\R$, uniquely determined $\RHO_1$-a.e.\ up to a constant, such that $\gamma=(\ID\times\nabla\phi) \#(\RHO_1\LEB^d)$; see \cites{GangboMcCann, Caffarelli, AmbrosioGigliSavare, Villani}. This identity implies that
$$
    \WAS(\RHO_1,\RHO_2)^2 = \int_{\R^d} |\nabla\phi(x)-x|^2 \RHO_1(x) \,dx
    \quad\text{and}\quad
    (\nabla\phi)\#(\RHO_1\LEB^d) = \RHO_2\LEB^d.
$$
We call $\nabla\phi$ an optimal transport map. It is
invertible $\RHO_1$-a.e.\ and the Hessian $D^2\phi$ exists and is
positive definite $\RHO_1$-a.e. When $d=1$, the Wasserstein distance
and the transport map can be computed more explicitly; see below.
\medskip

(2) We introduce a differentiable structure on $\SP(\R^d)$ as follows:
For any point $\RHO\in \SP(\R^d)$, the tangent space
$\TAN_\RHO\SP(\R^d)$ is defined as the closure of the space of smooth
gradient vector fields in the $\L^2(\R^d, \RHO)$-norm. This definition
is motivated by the fact that for any absolutely continuous curve
$t\mapsto\RHO_t\in\SP(\R^d)$ with $\RHO_0=\RHO$, there exists a unique
$\BU\in\TAN_\RHO\SP(\R^d)$ with the property that
\begin{equation}
    \partial_t \RHO_t \big|_{t=0} + \nabla\cdot(\RHO\BU) = 0
    \quad\text{in $\D'(\R^d)$.}
\label{E:TRA}
\end{equation}
That is, for all test functions $\phi\in\D(\R^d)$ we have
$$
    \frac{d}{dt} \bigg|_{t=0} \int_{\R^d} \RHO_t(x) \phi(x) \,dx
        = \int_{\R^d} \BU(x)\cdot\nabla\phi(x) \RHO(x) \,dx,
$$
and so a change in density can be accounted for by a flux of mass along
the velocity field  $\BU$. This structure makes $\SP(\R^d)$ a
Riemannian manifold, and we define
$$
    \TAN\SP(\R^d) := \bigcup\Big\{ (\RHO,\BU)\colon \RHO\in\SP(\R^d),
        \BU\in\TAN_\RHO\SP(\R^d) \Big\}.
$$

(3) If $\BU = -\RHO^{-1}\nabla P(\RHO)$, then \eqref{E:TRA}
yields the porous medium equation \eqref{E:POR} at one instant in time.
This vector field is the ``gradient'' of the internal energy
\begin{equation}
    \INT[\RHO] := \int_{\R^d} U(\RHO(x))\,dx
    \quad\text{with}\quad
    P(\RHO) = U'(\RHO)\RHO-U(\RHO)
\label{E:INTERNAL}
\end{equation}
in the sense that $\BU$ is the uniquely determined element of minimal
length in the {\em sub\-differential} of $\INT[\RHO]$ with respect to
the Wasserstein distance. The vector field $\BU$ is indeed tangent
to $\SP(\R^d)$ because $-\RHO^{-1}\nabla P(\RHO) = -\nabla U'(\RHO)$ is
a gradient.

\medskip

The internal energies for the porous medium equation that arise most
frequently in physical applications are identical to those of
the isentropic Euler equations, namely \eqref{E:POLYT} or \eqref{E:ISO}. For the latter choice with $\kappa=1$, equation \eqref{E:POR} becomes
\begin{equation}
    \partial_t\RHO - \Delta\RHO = 0
    \quad\text{in $[0,\infty)\times\R^d$,}
\label{E:HE}
\end{equation}
i.e., the porous medium equation reduces to the heat equation.

The interpretation of \eqref{E:POR} as an abstract gradient flow
suggests a natural time discretization: Given a time step $\tau>0$ and
the value $\RHO^n\in\SP(\R^d)$ of the approximate solution at time
$t^n:=n\tau$, the value at time $t^{n+1}$ is defined as
\begin{equation}
    \RHO^{n+1} := \ARGMIN\Bigg\{
        \frac{1}{2\tau} \WAS(\RHO^n,\RHO)^2 + \INT[\RHO]
            \colon \RHO\in\SP(\R^d) \Bigg\}.
\label{E:DISGF}
\end{equation}
We assume that $\RHO^0$ has finite internal energy $\INT[\RHO^0] <
\infty$. It can then be shown that for all $n\GS 0$, problem
\eqref{E:DISGF} admits a unique solution with $\INT[\RHO^{n+1}] \LS \INT[\RHO^n]$. We describe the motivation for the minimization \eqref{E:DISGF} in Appendix~\ref{S:EU:LA}.

It is useful to reformulate the porous medium equation \eqref{E:POR} in Lagrangian terms. We fix a reference density $\bar{\RHO}$ and consider a time-dependent family of transport maps $X(t,\cdot) \colon \R^d \longrightarrow\R^d$ such that $X(t,\cdot)\#\bar{\RHO} = \RHO(t,\cdot)\LEB^d$ for all $t\GS 0$. If
$$
    X(t,x)=X(t,y)
    \quad\text{implies that}\quad
     \partial^+_t X(t,x)=\partial^+_t X(t,y)
$$
for $\bar{\RHO}$-a.e.\ $x,y\in\R^d$, where $\partial^+_t$ denotes the
forward partial derivative in time, then we can define an Eulerian
velocity $\BU$ by the formula $\partial^+_t X =: \BU\circ X$. For the
porous medium flow, the velocity is given by Darcy's law $\BU =
-\nabla U'(\RHO)$. This vector field is related to the first variation
of the functional $\mc{V}[X] := \INT[X\# \bar{\RHO}]$, which is
defined for all $\bar{\RHO}$-measurable maps
$X\colon\R^d\longrightarrow\R^d$.  Let $\frac{\delta
  \mc{V}[X]}{\delta X} := \nabla U'(\RHO)\circ X$. As explained in
Appendix~\ref{S:EU:LA}, we then have
\begin{align*}
    \frac{d}{d\EPS}\bigg|_{\EPS=0} \mc{V}[X+\EPS (\zeta\circ X)]
        &= \int_{\R^d} \frac{\delta\mc{V}[X]}{\delta X} \cdot
            (\zeta \circ X) \,d\bar{\RHO}
        = \int_{\R^d} \nabla U'(\RHO) \cdot \zeta \,\RHO \,dx
\end{align*}
for all smooth test functions $\zeta\colon\R^d\longrightarrow\R^d$. Recall that $X\#\bar{\RHO} = \RHO\LEB^d$.

The minimization \eqref{E:DISGF} can now be rewritten in terms of transport maps as
\begin{equation}
  X^{n+1} := \ARGMIN_X \Bigg\{
    \frac{1}{2\tau} \int_{\R^d} |X-X^n|^2 \,d\bar{\RHO}
        + \mc{V}[X] \Bigg\},
\label{E:DISGF3}
\end{equation}
where $X^0$ is chosen such that $X^0\#\bar{\RHO} = \RHO^0\LEB^d$.
We obtain $X^{n+1} = \BR^{n+1}\circ X^n$ for all $n\GS 0$, where $\BR^{n+1}$ is the optimal transport map pushing $\RHO^n\LEB^d$ forward to $\RHO^{n+1}\LEB^d$. The discrete analogue of Darcy's law is then $\BU^{n+1} = -\nabla U'(\RHO^{n+1})$, where
\begin{equation}
    \BU^{n+1} := \frac{\BR^{n+1}-\ID}{\tau} \circ (\BR^{n+1})^{-1}
        = \frac{X^{n+1}-X^n}{\tau} \circ (X^{n+1})^{-1}.
\label{E:UN1:DEF}
\end{equation}
We refer the reader to Appendix~\ref{S:EU:LA} for details. In particular, we have that
$$
    \frac{X^{n+1}-X^n}{\tau} = -\nabla U'(\RHO^{n+1})\circ X^{n+1}
        = -\frac{\delta\mc{V}[X]}{\delta X}\bigg|_{X=X^{n+1}},
$$
which shows that the minimization problem \eqref{E:DISGF3} is equivalent to one step of the backward Euler method for the dynamical system $\partial_t X = -\frac{\delta\mc{V}[X]}{\delta X}$.

Notice that if for each timestep we {\em reset the reference measure} $\bar{\RHO}$ to be equal to $\RHO^n\LEB^d$, then we may choose $X^n=\ID$ and obtain $X^{n+1} =\BR^{n+1}$, so the minimization \eqref{E:DISGF3} reduces to \eqref{E:DISGF} if the latter is rewritten as a minimization problem for the optimal transport map pushing $\RHO^n\LEB^d$ forward to $\RHO^{n+1}\LEB^d$; see again Appendix~\ref{S:EU:LA}. This is what we do in our first order schemes detailed below. In designing more complicated multi-step or Runge-Kutta schemes, however, it is useful to describe all the intermediate steps or stages using one single reference density, which requires the generality of the second half of \eqref{E:UN1:DEF}. We also emphasize that in contrast to the space of probability measures $\SP(\R^d)$, the space of transport maps $X\colon\R^d \longrightarrow\R^d$ is {\em flat}, with a linear structure and a translation invariant metric.
\medskip

A similar steepest descent interpretation as for the porous medium
equation can be given for the system of isentropic Euler equations
\eqref{E:IEE}. The philosophy here is that among all possible weak
solutions of \eqref{E:IEE} we try to pick the one that dissipates the
total energy \eqref{E:TOTL} as fast as possible, as advocated by
Dafermos \cite{Dafermos}. Note that unlike in the porous medium case,
for \eqref{E:IEE} the energy dissipation will be singular in the sense
that it only occurs when the solution is discontinuous.

At the heart of the variational time discretization introduced in
\cite{GangboWestdickenberg} lies an optimization problem similar to
\eqref{E:DISGF}, but with a different homogeneity in the timestep; see
\eqref{E:INTRELAX} below. In fact, the relaxation of internal energy is
a second order effect only. What dominates the flow is the transport of
mass along the integral curves of the velocity field. Let us recall the
time discretization introduced in \cite{GangboWestdickenberg}.

\begin{definition}[Variational Time Discretization]\label{D:VTD}
\mbox{}
\medskip

(0) Let initial data $(\RHO^0,\BU^0) \in \TAN\SP(\R^d)$ and $\delta>0$
be given. Put $t^0 := 0$.
\medskip

\noindent Assume that $\RHO^n\in\SP(\R^d)$ and $\BU^n \in
\L^2(\R^d,\RHO^n)$ are given, where $n\GS 0$. Then
\medskip

(1) Choose $\tau\in[\delta/2,\delta]$ in such a way that the
push-forward measure
$$
    (\ID+\tau\BU^n) \# (\RHO^n\LEB^d) =: \hat{\RHO}^n \LEB^d
$$
is absolutely continuous with respect to the Lebesgue measure. Find
\begin{equation}
    \hat{\BU}^n \in \L^2(\R^d,\RHO^n)
    \quad\text{such that}\quad
    \left\{\begin{gathered}
        \vphantom{\int}
            (\ID+\tau\hat{\BU}^n)\#(\RHO^n\LEB^d)
                = \hat{\RHO}^n\LEB^d
\\
        \int_{\R^d} |\tau\hat{\BU}^n|^2 \RHO^n \,dx
            = \WAS(\RHO^n,\hat{\RHO}^n)^2
    \end{gathered}\right\}.
\label{E:OTU}
\end{equation}
That is, find the velocity field $\hat{\BU}^n$ with minimal $\L^2(\R^d,
\RHO^n)$-norm such that the map $\ID+\tau\hat{\BU}^n$ pushes the
measure $\RHO^n \LEB^d$ forward to $\hat{\RHO}^n\LEB^d$. Such a vector
field always exists and is uniquely determined. It is given by the
gradient of a semi-convex function, and is therefore a tangent vector.
Moreover, we have the estimate
\begin{equation} \label{E:MAX:DISSIPATE}
    \int_{\R^d} |\hat{\BU}^n|^2 \RHO^n \,dx
        \LS \int_{\R^d} |\BU|^2 \RHO^n \,dx
\end{equation}
for any velocity $\BU$ satisfying the constraint that $(\ID+\tau\BU)\#
\RHO^n = \hat{\RHO}^n$.\ This shows that by replacing
$\BU^n$ by the optimal transport velocity $\hat{\BU}^n$, we decrease
the kinetic energy as much as possible. The map $\ID+\tau\hat{\BU}^n$
is invertible $\hat{\RHO}^n$-a.e.
\medskip

(2) Update the density by computing the minimizer
\begin{equation}
    \RHO^{n+1} := \ARGMIN\Bigg\{
        \frac{3}{4\tau^2} \WAS(\hat{\RHO}^n,\RHO)^2 + \INT[\RHO]
            \colon \RHO\in\SP(\R^d) \Bigg\}.
\label{E:INTRELAX}
\end{equation}
Here $\INT[\RHO]$ is the internal energy functional from \eqref{E:INTERNAL}, which is a part of \eqref{E:TOTL}. One can show that the density $\RHO^{n+1}$ is uniquely determined and that
\begin{gather*}
    \bigg(\ID+\frac{2\tau^2}{3}\nabla U'(\RHO^{n+1})\bigg)
            \#\RHO^{n+1}
        = \hat{\RHO}^n,
\\[1ex]
    \int_{\R^d} \bigg|\frac{2\tau^2}{3}\nabla U'(\RHO^{n+1})\bigg|^2
            \RHO^{n+1} \,dx
        = \WAS(\RHO^{n+1},\hat{\RHO}^n)^2.
\end{gather*}
The density $\RHO^{n+1}$ is thus regular in the sense that $\frac{2
\tau^2}{3} \nabla U'(\RHO^{n+1}) \in \L^2(\R^d, \RHO^{n+1})$. In fact,
the latter vector field is the gradient of a semi-convex function.
\medskip

(3) Update the velocity by defining
\begin{equation}
     \BU^{n+1} := \hat{\BU}^n \circ (\ID+\tau\hat{\BU}^n)^{-1}
        \circ \bigg(\ID+\frac{2\tau^2}{3} \nabla U'(\RHO^{n+1})\bigg)
            -\tau\nabla U'(\RHO^{n+1}).
\label{E:UN:UP}
\end{equation}

(4) Let $t^{n+1} := t^n + \tau$, increase $n$ by one, and continue with
Step~(1).
\end{definition}

We refer the reader to \cite{GangboWestdickenberg} for more details. It
is shown there that
$$
    E[\RHO^{n+1},\BU^{n+1}] \LS E[\RHO^n,\BU^n]
    \quad\text{for all $n\GS 0$.}
$$
One can then define the function
$$
    (\RHO,\BU)(t, \cdot) := (\RHO^n,\BU^n)
    \quad\text{for all $t\in[t^n,t^{n+1})$ and $n\GS 0$,}
$$
which is piecewise constant in time and approximates a weak,
energy-dissipating solution of the isentropic Euler equation
\eqref{E:IEE}.\ A formal argument is given in
\cite{GangboWestdickenberg} for the convergence of this approximation
when the timestep $\delta\rightarrow 0$.
\medskip

We give a derivation of the time discretization of Definition~\ref{D:VTD} in Appendix~\ref{S:DERIVE:VTD}. From the point of view of numerical
analysis, it is interesting to generalize this scheme to a family of
first order implicit methods that also includes the backward Euler
method as a special case.  To this end, let us temporarily assume the
solution is smooth and write the isentropic Euler equations in a
Lagrangian frame as
\begin{equation}
    \frac{d}{dt}
    \begin{pmatrix}
        X \\
        V
    \end{pmatrix}
    =
    \begin{pmatrix}
        V \\
        \mb{f}[X]
    \end{pmatrix}
    \quad\text{with}\quad
    \mb{f}[X]
        := -\frac{\delta\mc{V}[X]}{\delta X}
        = -\nabla U'(\RHO_X)\circ X,
\label{E:LAGRANGIAN:FORMULATION}
\end{equation}
where $X\#\bar{\RHO} =: \RHO_X\LEB^d$ and $V := \BU\circ X$ is the pull-back of the Eulerian velocity. Again we set $\mc{V}[X] := \INT[X\#\bar{\RHO}]$, as in the porous medium case discussed above. Let now $\alpha\in(0,1]$ be a fixed parameter and consider the implicit method
\begin{equation}
    \begin{pmatrix}
        X^{n+1} \\
        V^{n+1}
    \end{pmatrix}
    =
    \begin{pmatrix}
        X^{n} \\
        V^{n}
    \end{pmatrix}
    + \tau
    \begin{pmatrix}
        (1-\alpha) V^{n} + \alpha V^{n+1} \\
        \mb{f}[X^{n+1}]
    \end{pmatrix}.
\label{E:FAMILY:OF:SCHEMES}
\end{equation}
Using the first equation to eliminate $V^{n+1}$ in the second,
we obtain
\begin{gather}
    \frac{1}{\alpha\tau^2} \Big( X^{n+1} - (X^n + \tau V^n) \Big)
        = \mb{f}[X^{n+1}],
\label{E:UPDATE:LA:A}\\
    V^{n+1} = V^n + \tau\mb{f}[X^{n+1}].
\label{E:UPDATE:LA:B}
\end{gather}
Equation \eqref{E:UPDATE:LA:A} will be satisfied by the solution of the
minimization problem
\begin{equation} \label{E:AMIN:LA}
    X^{n+1} = \ARGMIN_X \Bigg\{
        \frac{1}{2\alpha\tau^2} \int_{\R^d} |X-(X^n+\tau V^n)|^2 \,d\bar{\RHO}
            + \mc{V}[X] \Bigg\}.
\end{equation}
Once this is solved, we substitute the left-hand side of
\eqref{E:UPDATE:LA:A} for $\mb{f}[X^{n+1}]$ in \eqref{E:UPDATE:LA:B} to obtain
$V^{n+1}$. This is the Lagrangian reformulation of the energy minimization in Step~(2) of Definition~\ref{D:VTD}. As explained there, we combine this step with a velocity projection, so we replace in each step the velocity $V^n$ by an optimal transport velocity, which in the Lagrangian framework takes the form $\hat{V}^n = \hat{\BU}^n\circ X^n$ with $\hat{\BU}^n$ given by \eqref{E:OTU}. This modification dissipates the maximum amount of kinetic energy
without changing the convected distribution of mass.

The three choices of $\alpha$ that seem most natural are $\alpha=1$,
which corresponds to the backward Euler method; $\alpha=\frac{2}{3}$,
which corresponds to the original time discretization proposed
in \cite{GangboWestdickenberg} (see Appendix~\ref{S:DERIVE:VTD} for further information); and $\alpha=\frac{1}{2}$, which causes
the acceleration $A^{n+1}=\mb{f}[X^{n+1}]$ to agree with its
Taylor expansion, $A^{n+1}\approx 2\tau^{-2}\big(X^{n+1} -
(X^n + \tau V^n)\big)$.  All three of these variants yield first
order methods that appear to capture shocks and rarefaction waves
correctly in our numerical experiments. We present a second order
version in Section~\ref{SS:SOS}.


\section{Variational Particle Scheme}\label{S:VPS}

We now introduce a fully discrete version of the variational time
discretization discussed in the previous section for the
one-dimensional isentropic Euler equations
\begin{equation}
\begin{aligned}
    & \partial_t \RHO + \partial_x(\RHO u) = 0
\\
    & \partial_t ( \RHO u ) + \partial_x ( \RHO u^2 + P(\RHO) )
        = 0
\end{aligned}
\quad\text{in $[0,\infty)\times\R$,}
\label{E:IE}
\end{equation}
as well as for the one-dimensional porous medium equation
\begin{equation}
    \partial_t \RHO - \Delta P(\RHO) = 0
    \quad\text{in $[0,\infty)\times\R$.}
\label{E:PMEQ}
\end{equation}
The pressure $P$ is related to the internal energy by equation
\eqref{E:PRESS}, with energy density $U$ given by either
\eqref{E:POLYT} in the case of polytropic gases (we put $\kappa=1$ for
the porous medium equation), or by \eqref{E:ISO} in the case of
isothermal gases.


\subsection{First Order Scheme}\label{SS:FOS}

In the simplest version of the algorithm, we assume the fluid to be
composed of $N$ particles of equal mass $m:=1/N$, where $N>1$. At a
given time $t^n$, the particles are located at positions $\BX^n =
\{x^n_1,\ldots,x^n_N\}$ with velocities $\BU^n = \{u^n_1,\ldots,
u^n_N\}$. In the following, we will always assume that the $x^n_i$ are
ordered and {\em strictly increasing}. The total energy is then given
by
$$
    E(\BX^n,\BU^n) := \sum_{i=1}^N \frac{1}{2}m |u^n_i|^2
        + \sum_{i=1}^{N-1} U\bigg(\frac{m}{x^n_{i+1}-x^n_i}\bigg)
            (x^n_{i+1}-x^n_i).
$$
In accordance with thermodynamics, we think of the density as the
inverse of the specific volume, which we define as the distance between
two neighboring particles divided by the particle mass. This formula is
used in the internal energy.

\begin{definition}[Variational Particle Scheme, Version 1]\label{D:SCHEME1}
\mbox{}
\medskip

(0) Let initial positions $\BX^0$, velocities $\BU^0$, and $\tau>0$ be
given. Put $t^0 := 0$.
\medskip

\noindent Assume that positions $\BX^n$ and velocities $\BU^n$ are
given for $n\GS 0$. Then
\medskip

(1) Compute intermediate positions defined by $y_i := x^n_i + \tau
u^n_i$ for all $1\LS i\LS N$. Find the permutation $\sigma$ of the set
$\{1,\ldots,N\}$ with the property that the $\hat{x}^n_i :=
y_{\sigma(i)}$ are nondecreasing in $i$. We assume that if there exist
indices $i<j$ with $\hat{x}_i = \hat{x}_j$, then $\sigma(i) <
\sigma(j)$. Then define $\hat{u}^n_i := (\hat{x}^n_i-x^n_i)/\tau$ for
all $i$.
\medskip

(2) Define new positions $\BX^{n+1}$ as the minimizer of the functional
\begin{equation}
    F(\BZ) := \sum_{i=1}^N \frac{3}{4\tau^2} m |z_i-\hat{x}^n_i|^2
        + \sum_{i=1}^{N-1} U\bigg(\frac{m}{z_{i+1}-z_i}\bigg)
            (z_{i+1}-z_i)
\label{E:MIINT}
\end{equation}
over all families of positions $\BZ\in\R^N$ such that $z_{i+1}>z_i$ for
all $i$.
\medskip

(3) Compute the new velocities $u^{n+1}_i := \hat{u}^n_i +
\frac{3}{2\tau}(x^{n+1}_i - \hat{x}^n_i)$ for all $i$.
\medskip

(4) Increase $n$ by one and continue with Step~(1).
\end{definition}

\begin{remark}
In Step~(1) we perform a free transport of all particles in the
direction of their given velocities. We obtain a new distribution of
mass that is characterized by the new particle positions $y_i$. Note
that we defined the density in terms of the distance between {\em
neighboring} particles. Since the original ordering can be destroyed
during the free transport step, we rearrange the new positions in
nondecreasing order. On the discrete level, it is not necessary to make
sure that positions do not coincide because we never compute the
internal energy of this intermediate configuration. Therefore we can
use the same timestep $\tau$ for all updates.

To understand the definition of the new velocities in Step~(1), notice
that the Wasserstein distance as defined in \eqref{E:WAS} can easily be
generalized to pairs of {\em probability measures} $\mu_1$ and $\mu_2$
with finite second moment; see \cites{AmbrosioGigliSavare, Villani} for
example. If the two measures are given as convex combinations of Dirac
measures:
$$
    \mu_1 = \frac{1}{N} \sum_{i=1}^N \delta_{x_i}
    \quad\text{and}\quad
    \mu_2 = \frac{1}{N} \sum_{i=1}^N \delta_{y_i},
$$
where $x_i\in\R^d$ and $y_i\in\R^d$ for all $i$, then the squared
Wasserstein distance between the measures $\mu_1$ and $\mu_2$ coincides
with the minimum of the functional
$$
    W(\sigma)^2 := \frac{1}{N} \sum_{i=1}^N |x_i-y_{\sigma(i)}|^2
$$
among all permutations $\sigma$ of the index set $\{1,\ldots,N\}$. This
amounts to a Linear Assignment Problem, which in the one-dimensional
case reduces to sorting the positions $y_i$. This can be done with
worst-case complexity $\mathcal{O}(N\log N)$. In our case, the $y_i$
are actually largely ordered already, so the sorting is quite
inexpensive. Step~(1) therefore amounts to computing the Wasserstein
distance between the old and the intermediate positions of particles,
and then finding the optimal velocity that achieves this transport: We
have that $x^n_i + \tau \hat{u}^n_i = \hat{x}^n_i$ for all $i$.

The functional $F$ in \eqref{E:MIINT} is the analogue of the one in
\eqref{E:INTRELAX}. Since we minimize over positions $z_i$ that are
increasing, the first term of \eqref{E:MIINT} can again be interpreted
as a Wasserstein distance squared between convex combinations of Dirac
measures located at positions $z_i$ and $\hat{x}^n_i$. Note that the
choice $z_i := x^n_i$ for all $i$ is admissible in \eqref{E:MIINT} and
results in a finite value of the functional. Therefore the new
positions $x^{n+1}_i$ must have finite internal energy and thus be
strictly increasing.
\end{remark}

\begin{remark}
With only minor modifications we obtain a particle scheme for the
porous medium equation \eqref{E:PMEQ}, which contains the heat equation
as a special case. Since now the only unknown is the density $\RHO$,
Steps~(1) and (3) in Definition~\ref{D:SCHEME1} become irrelevant. In
Step~(2), we substitute the factor $1/(2\tau)$ for $3/(4\tau^2)$
to match \eqref{E:DISGF}, and the intermediate positions $\hat{x}_i$
must be replaced by $x_i$.
Of course, we can also implement the backward Euler method for the
isentropic and isothermal Euler equations by changing $3/(4\tau^2)$ and
$3/(2\tau)$ to $1/(2\tau^2)$ and $1/\tau$, respectively.
\end{remark}


\subsection{Second Order Scheme}\label{SS:SOS}

In this section, we derive a variant of the scheme of
Definition~\ref{D:SCHEME1} that is formally second order in space and
time. Before doing so, let us collect a few facts about the Wasserstein
distance in one space dimension. For any probability measure $\mu$ in
$\R$ we define the distribution function
$$
    F_\mu(t) := \mu\big((-\infty,t))
    \quad\text{for all $t\in\R$.}
$$
The function $F_\mu$ is nondecreasing, but not necessarily strictly,
and can be discontinuous if the measure $\mu$ contains atoms (Dirac
measures). We define
$$
    F_\mu^{-1}(s) := \sup\Big\{ t\in\R \colon F_\mu(r) \LS s \Big\}
    \quad\text{for all $s\in[0,1]$,}
$$
which is the generalized inverse of $F_\mu$. For any pair of
probability measures $(\mu,\nu)$ on $\R$ with finite second moment, the
Wasserstein distance between $\mu$ and $\nu$ can then be computed in
terms of the inverse distribution functions as
\begin{equation}
    \WAS(\mu,\nu)^2
        := \int_{[0,1]} |F_\mu^{-1}(s)-F_\nu^{-1}(s)|^2 \,ds.
\label{E:ONEDWAS}
\end{equation}
If $\mu$ does not contain any atoms, then $F_\mu$ is continuous, and
the composition $F_\nu^{-1}\circ F_\mu$ is an optimal transport map
pushing $\mu$ forward to $\nu$; see Theorem~6.0.2 in
\cite{AmbrosioGigliSavare}. Note that the inverse function $F_\mu^{-1}$
pushes the measure $\IND_{[0,1]}\LEB^1$ forward to $\mu$.

Assume now that the measure $\mu$ is absolutely continuous with respect
to the Lebesgue measure, with a piecewise constant density: Let numbers
$x_0 < \ldots < x_N$ and $m_i\GS 0$ be given such that $\sum_{i=1}^N
m_i = 1$. Let $\mu = \RHO\LEB^1$ with
\begin{equation}
    \RHO(x) := \begin{cases}
        0 & \text{if $x<x_0$ or $x\GS x_N$,}
\\
        \DST\frac{m_i}{x_i-x_{i-1}}
            & \text{if $x\in [x_{i-1},x_i)$ and
                $1\LS i\LS N$.}
    \end{cases}
\label{E:DEFRHO}
\end{equation}
Then the inverse distribution function of $\mu$ is a piecewise linear
function that can be computed explicitly. Defining $s_k := \sum_{i=1}^k
m_i$ for $0\LS k\LS N$, we have
$$
    F_\mu^{-1}(s)
        = x_{k-1} \frac{s_k-s}{s_k-s_{k-1}}
        + x_k \frac{s-s_{k-1}}{s_k-s_{k-1}}
    \quad\text{for all $s\in[s_{k-1},s_k]$}
$$
and $1\LS k\LS N$. If the second measure $\nu$ is piecewise constant as
well, it is possible to compute the Wasserstein distance between $\mu$
and $\nu$ using formula \eqref{E:ONEDWAS}. This is particularly simple
if both $\mu$ and $\nu$ have the same number of intervals, and if the
same mass is assigned to each interval, because then the inverse
distribution functions $F_\mu^{-1}$ and $F_\nu^{-1}$ are piecewise
linear on the same intervals.

Assume now that $\nu=\hat{\RHO}\LEB^1$ with
$$
    \hat{\RHO}(x) := \begin{cases}
        0 & \text{if $x<\hat{x}_0$ or $x\GS \hat{x}_N$,}
\\
        \DST\frac{\hat{m}_i}{\hat{x}_i-\hat{x}_{i-1}}
            & \text{if $x\in [\hat{x}_{i-1},\hat{x}_i)$ and
                $1\LS i\LS N$,}
    \end{cases}
$$
for suitable numbers $\hat{x}_0 < \ldots < \hat{x}_N$ and $\hat{m}_i\GS
0$ with $\sum_{i=0}^N \hat{m}_i = 1$.\ We want to project $\nu$ onto
the space of measures of the form $\mu=\RHO\LEB^d$, where the density
$\RHO$ is given by \eqref{E:DEFRHO} and the masses $m_i$ are fixed.
More precisely, we want to choose the positions $\BX := (x_0, \ldots,
x_N)$ in \eqref{E:DEFRHO} so that the Wasserstein distance
$$
    \WAS(\mu,\nu)^2 = \int_{[0,1]} |F_\mu^{-1}(s)-F_\nu^{-1}(s)|^2 \,ds
$$
is minimal. This amounts to a quadratic minimization problem that can
be solved easily. Let $\{\varphi_k\}_{k=0}^N$ be the standard finite
element hat functions with vertices $s_k$ as defined above. Then the
minimizer is given by $\BX = A^{-1}b$, where
$$
    A_{kl} := \int_{[0,1]} \varphi_k(s)\varphi_l(s) \,ds
    \quad\text{and}\quad
    b_k := \int_{[0,1]} \varphi_k(s) F_\nu^{-1}(s) \,ds
$$
for all $0\LS k,l\LS N$. The calculation of $A_{kl}$ and $b_k$ is
straightforward. Notice that the function $F_\nu^{-1}$ is piecewise
linear on the intervals $[\hat{s}_{k-1},\hat{s}_k]$ with $\hat{s}_k :=
\sum_{i=1}^k \hat{m}_i$, and the $\hat{s}_k$ are different from $s_k$
(otherwise there is nothing to do). The computation of $b_k$ therefore
requires a partition of $[0,1]$ that uses both sets of positions. Even
if the $\hat{x}_i$ are nondecreasing, the components $x_i$ of the
minimizer may not be.

The projection outlined in the previous paragraph works also in the
degenerate situation where $\hat{x}_{i-1} = \hat{x}_i$ for one or more
$i$. In this case $\nu$ has a Dirac measure at position $\hat{x}_i$ and
$\hat{m}_i$ denotes the mass of that measure. The inverse distribution
function $F_\nu^{-1}$ is then no longer continuous, but it is still
straightforward to compute the numbers $b_i$. The minimizer $\BX$ can
be obtained as above.
\medskip

We can now define our second order method for the isentropic Euler
equations. Let us first discuss the accuracy in space. Instead of
working with point masses, we now approximate the density by a function
that is piecewise constant on intervals of variable length.
Specifically, we fix numbers $m_i\GS 0$ with $\sum_{i=1}^N m_i = 1$. At
any time $t^n$, the density $\RHO^n$ is then determined by a vector
$\BX^n = (x^n_0,x^n_1,\ldots, x^n_N) \in \R^{N+1}$ of positions with
$x^n_{i-1}<x^n_i$ for all $1\LS i\LS N$, through the formula
\begin{equation}
    \RHO^n(x) := \begin{cases}
        0 & \text{if $x<x^n_0$ or $x\GS x^n_N$,}
\\
        \DST\frac{m_i}{x^n_i-x^n_{i-1}}
            & \text{if $x\in [x^n_{i-1},x^n_i)$ and
                $1\LS i\LS N$.}
    \end{cases}
\label{E:RHO:N}
\end{equation}
We approximate the velocity at time $t^n$ by a piecewise linear
function on the intervals $[x^n_{i-1}, x^n_i)$. It is determined by
$\BU^n = (u^n_0,u^n_1,\ldots,u^n_N) \in \R^{N+1}$ as
$$
    u^n(x) := \begin{cases}
        0 & \text{if $x<x^n_0$ or $x\GS x^n_N$,}
\\
        \DST\frac{x^n_i-x}{x^n_i-x^n_{i-1}} u^n_{i-1}
                + \frac{x-x^n_{i-1}}{x^n_i-x^n_{i-1}}u^n_i
            & \text{if $x\in [x^n_{i-1},x^n_i)$ and
                $1\LS i\LS N$.}
    \end{cases}
$$

For any $n\GS 0$ consider now the push-forward measure $\nu :=
(\ID+\tau u^n)\# (\RHO^n\LEB^1)$. Let $y_i := x^n_i + \tau u^n_i$ for
all $0\LS i\LS N$, and let $\sigma$ be a permutation of $\{0,1,\ldots,
N\}$ such that the positions $\hat{x}_i := y_{\sigma(i)}$ are
nondecreasing in $i$. We assume that if there are indices $i < j$ with
$\hat{x}_i = \hat{x}_j$, then $\sigma(i) < \sigma(j)$. The measure
$\nu$ is piecewise constant on intervals $[\hat{x}_{i-1}, \hat{x}_i)$
for which $\hat{x}_{i-1} < \hat{x}_i$, and we denote the mass carried
by this interval by $\hat{m}_i$ and the density by
$\hat{\RHO}_i := \hat{m}_i / (\hat{x}_i-\hat{x}_{i-1})$.
If $\hat{x}_{i-1} = \hat{x}_i$,
then the measure $\nu$ has a Dirac measure located at position
$\hat{x}_i$ and $\hat{m}_i$ denotes the mass of this measure. By
construction, we have $\hat{m}_i\GS 0$ and $\sum_{i=1}^N \hat{m}_i =
1$.

To compute the numbers $\hat{m}_i$ we first initialize them to zero. We
then run through the original intervals $[x_{i-1},x_i)$ and assign
their masses to the new intervals they are mapped onto, in proportion
to the lengths of these new intervals:
\medskip
\begin{tabbing}
\hspace*{2em} \= \hspace*{2em} \= \hspace*{2em} \= \hspace*{2em} \= \kill
    \> \textbf{for} $i=1$ \textbf{to} $N$
\\
    \>\> $k = \min\{\sigma^{-1}(i-1),\sigma^{-1}(i)\}$
\\
    \>\> $l = \max\{\sigma^{-1}(i-1),\sigma^{-1}(i)\}$
\\
    \>\> \textbf{for} $j=k+1$ \;\textbf{to}\; $l$
\\
    \>\>\> \textbf{if} $\hat{x}^n_k = \hat{x}^n_l$
\\
    \>\>\>\> $\hat{m}_j = \hat{m}_j + \DST\frac{m_i}{l-k}$
        \quad (point mass case)
\\
    \>\>\> \textbf{else}
\\
    \>\>\>\> $\hat{m}_j = \hat{m}_j +
        \DST m_i\frac{\hat{x}^n_j-\hat{x}^n_{j-1}}
            {\hat{x}^n_l - \hat{x}^n_k}$
    \quad (distribute according to length).
\end{tabbing}
\medskip

As we explained above, the computation of the Wasserstein distance in
one space dimension is very simple if the pair of measures is induced
by two piecewise constant densities that have the same number $N$ of
intervals, and that assign the same mass $m_i$ to the $i$th interval.
We therefore project the push-forward $\nu = (\ID+\tau u^n)\#
(\RHO^n\LEB^1)$ onto the space of measures of this form, before
proceeding with the variational time discretization (VTD) algorithm of
Definition~\ref{D:VTD}. The minimization in \eqref{E:DISGF} is only
performed over densities $\RHO$ of the same form; see below for more
details. Otherwise this version of the fully discrete algorithm is
identical to the VTD algorithm.
\medskip

Next we show how to modify the algorithm to achieve second order
accuracy in time. Consider again the Lagrangian formulation
\eqref{E:LAGRANGIAN:FORMULATION} of the system of isentropic Euler equations.
The key feature of the scheme \eqref{E:FAMILY:OF:SCHEMES} that allows
$X^{n+1}$ to be determined by solving a minimization problem of the
form \eqref{E:AMIN:LA} is that the functional $\mb{f}$ is only evaluated at
$X^{n+1}$ in \eqref{E:FAMILY:OF:SCHEMES}.  It is also desirable to
choose an $L$-stable method to prevent oscillations from growing near
shocks.  The second order backward differentiation formula method
\cite{HairerNorsettWanner} satisfies both of these conditions. Therefore it is a natural candidate to try to adapt to our variational particle scheme framework.  Applied to problem \eqref{E:LAGRANGIAN:FORMULATION}, BDF2 takes the form
\begin{equation}
    \frac{3}{2}
    \begin{pmatrix}
        X^{n+1} \\
        V^{n+1}
    \end{pmatrix}
    -2
    \begin{pmatrix}
        X^n \\
        V^n
    \end{pmatrix}
    + \frac{1}{2}
    \begin{pmatrix}
        X^{n-1} \\
        V^{n-1}
    \end{pmatrix}
    =\tau
    \begin{pmatrix}
        V^{n+1} \\
        \mb{f}[X^{n+1}]
    \end{pmatrix}.
\label{E:BDF2:DEF}
\end{equation}
Solving for $X:=X^{n+1}$ and $V^{n+1}$ in turn, we obtain
\begin{gather}
  \frac{3}{\tau^2}\left[
        \Big( X - (X^n+\tau V^n) \Big)
            -\frac{1}{4}\Big( X - \big(X^{n-1}
                + 2\tau ( \TST\frac{2}{3}V^n +
                          \frac{1}{3}V^{n-1}) \big) \Big) \right]
    = \mathbf{f}[X],
\nonumber\\
    \vphantom{\bigg|}
    V^{n+1}
        = 2V^n - ( \TST\frac{2}{3}V^n + \frac{1}{3}V^{n-1} )
        +\frac{2\tau}{3}\mb{f}[X^{n+1}].
\label{E:BDF2:UPDATE}
\end{gather}
The first equation characterizes the solution of the minimization problem \eqref{E:F2} below, which is similar to \eqref{E:AMIN:LA}. As before, we
handle the formation of shocks by dissipating the maximum amount of
energy possible without changing the convected distribution of mass
from $t^{n-1}$ to $t^n$ to $t^{n+1}$. Notice that the convex
combination $\frac{2}{3}V^n+\frac{1}{3}V^{n-1}$ (rather than
$V^{n-1}$) is used to transport mass from time $t^{n-1}$ to time
$t^{n+1}$.  As $V$ and $\BU$ are related via $V=\BU\circ X$,
we will use $\BU$ to denote velocities below.
Here is our second order method for the isentropic Euler
equations \eqref{E:IEE}:

\begin{definition}[Variational Particle Scheme, Version 2]\label{D:SCHEME2}
\mbox{}
\medskip

\noindent Fix once and for all masses $m_i\GS 0$ with $\sum_{i=1}^N m_i
= 1$.
\medskip

(0) Let initial positions $\BX^0$, velocities $\BU^0$, and $\tau>0$ be
given. Put $t^0 := 0$.
\medskip

\noindent Assume that positions $\BX^n$ and velocities $\BU^n$ are
given for $n\GS 0$. Then
\medskip

(1a) Compute the intermediate positions $y_i := x^n_i + \tau u^n_i$ for
all $0\LS i\LS N$. Sort them and redistribute the masses to obtain
$\hat{x}^n_i$ and $\hat{m}_i$. Project back onto the space of densities
that are constant on $N$ intervals and assign mass $m_i$ to the $i$th
interval. Let $\BX = (x'_0, \ldots, x'_N)$ denote the positions of the
minimizer obtained from the quadratic minimization problem induced by
this projection. See the beginning of this section for more details.
Then define the velocities $u'_i := (x'_i-x^n_i)/\tau$.
\medskip

(1b) Compute the intermediate positions $y_i := x^{n-1}_i + 2\tau
(\frac{2}{3}u^n_i + \frac{1}{3}u^{n-1}_i)$ for all $0\LS i\LS N$. Sort
them and redistribute the masses to obtain $\hat{x}^n_i$ and
$\hat{m}_i$ as before. Project back onto the space of densities that
are constant on $N$ intervals and assign mass $m_i$ to the $i$th
interval. Let $\BX'' = (x''_0, \ldots, x''_N)$ denote the positions of
the minimizer obtained from the quadratic minimization problem induced
by this projection. Then define the velocities $u''_i :=
(x''_i-x^{n-1}_i)/(2\tau)$, which play the role of
$\frac{2}{3}V^n+\frac{1}{3}V^{n-1}$ in (\ref{E:BDF2:UPDATE}).

{\em Note: Skip Step~(1b) on the first iteration.}
\medskip

(2) Minimize the functional
\begin{equation}
    F(\BZ) := \frac{3}{2\tau^2} \|\BZ - \BX'\|^2_m
        -\frac{3}{8\tau^2} \|\BZ - \BX''\|^2_m
        + \sum_{i=1}^{N} U\bigg(\frac{m_i}{z_{i}-z_{i-1}}\bigg)
            (z_{i}-z_{i-1})
\label{E:F2}
\end{equation}
over all $\BZ = (z_0, \ldots, z_N)$ such that $z_{i-1}<z_i$ for all
$1\LS i\LS N$. We define
$$
    \|\BZ\|_m^2
        := \int_0^1 \Bigg| \sum_{i=0}^N z_i\varphi_i(s)
            \Bigg|^2 \,ds,
$$
where $\varphi_i$ is the piecewise linear function satisfying
$\varphi_i(s_k) = \delta_{ik}$ for all $0\LS k\LS N$, with $s_k :=
\sum_{i=1}^k m_i$. In fact, we have  $\|\BZ\|_m^2 = \BZ^T A\BZ$ with
the matrix $A$ defined by
$$
    A_{ij} :=\int_{[0,1]} \varphi_i(s) \varphi_j(s) \,ds
$$
for all $0\LS i,j\LS N$. We denote by $\BX^{n+1} = (x^{n+1}_0, \ldots,
x^{n+1}_N)$ the minimizer of this convex optimization problem. Notice
that the Euler-Lagrange equations of \eqref{E:F2} are of the form
\eqref{E:BDF2:UPDATE} in the second order backwards differentiation
scheme.

{\em Note: On the first iteration, minimize instead the functional}
$$
    F(\BZ) := \frac{3}{4\tau^2} \big\|\BZ-\BX'\big\|^2_m
        + \sum_{i=1}^{N} U\bigg(\frac{m_i}{z_{i}-z_{i-1}}\bigg)
            (z_{i}-z_{i-1}).
$$
%

(3) Compute the new velocities
$$
    u^{n+1}_i := 2u_i' - u_i''
        + \frac{2}{\tau} \Big( x^{n+1}_i-x'_i \Big)
        - \frac{1}{2\tau} \Big( x^{n+1}_i-x''_i \Big).
$$
for all $0\LS i\LS N$. Notice that this formula is of the form
\eqref{E:BDF2:UPDATE}, where we recall that $u_i''$ represents
$\frac{2}{3}u_i^n + \frac{1}{3}u_i^{n-1}$, not $u_i^{n-1}$.

{\em Note: On the first iteration, define instead}
$$
    u^{n+1}_i := u'_i
        + \frac{3}{2\tau}\Big( x^{n+1}_i - x'_i \Big).
$$
%

(4) Increase $n$ by one and continue with Step~(1).
\end{definition}

\begin{remark}
While the $\hat{x}^n_j$ are nondecreasing by construction, the
projection can cause the components of the minimizer to be slightly out
of order. This is acceptable as the minimization \eqref{E:INTRELAX}
over densities $\RHO$ of the form \eqref{E:RHO:N} gives the same result
if we replace the push-forward measure $\nu:=(\ID+\tau u^n)\# (\RHO^n
\LEB^1)$ by its projection onto the space of densities of the form
\eqref{E:RHO:N}. Thus, in spite of the projection, we still minimize
the distance from $\nu$ to compute $\RHO^{n+1}$. The values $x_i^{n+1}$
will again be strictly increasing since otherwise the internal energy
would be infinite.
\end{remark}

\begin{remark}
As before, the same procedure works for the porous medium equation if
we eliminate Steps~(1) and (3) and replace the functional $F$ in step
(2) by
$$
    F(\BZ) := \frac{1}{\tau} \|\BZ - \BX^n\|^2_m
        -\frac{1}{4\tau} \|\BZ - \BX^{n-1}\|^2_m
            + \sum_{i=1}^{N} U\bigg(\frac{m_i}{z_{i}-z_{i-1}}\bigg)
                (z_{i}-z_{i-1}).
$$
On the first iteration, one should use
$$
    F(\BZ) := \frac{1}{2\tau} \|\BZ - \BX^n\|^2_m
            + \sum_{i=1}^{N} U\bigg(\frac{m_i}{z_{i}-z_{i-1}}\bigg)
                (z_{i}-z_{i-1}).
$$
\end{remark}

\begin{remark}
We have also constructed second (and higher) order variational
particle schemes in the diagonally implicit Runge-Kutta framework
\cite{HairerNorsettWanner}. The underlying Runge-Kutta method must
be ``stiffly accurate'', i.e., the numerical solution at $t^{n+1}$
must agree with the last internal stage of the Runge-Kutta
procedure.  We illustrate the basic idea using the two-stage second
order DIRK scheme with Butcher array
\begin{equation}
    \begin{array}{c|cc}
        1/4 & 1/4 \\
        1 & 2/3 & 1/3 \\
        \hline & 2/3 & 1/3
    \end{array}.
\label{E:BUTCHER}
\end{equation}
The system we need to solve is
\begin{align}
    \begin{pmatrix}
        X^{n+\frac{1}{4}} \\
        V^{n+\frac{1}{4}}
    \end{pmatrix}
    &=
    \begin{pmatrix}
        X^{n} \\
        V^{n}
    \end{pmatrix}
    + \frac{\tau}{4}
    \begin{pmatrix}
        V^{n+\frac{1}{4}} \\
        \mb{f}[X^{n+\frac{1}{4}}]
    \end{pmatrix},
\label{E:STAGE1}\\
    \begin{pmatrix}
        X^{n+1} \\
        V^{n+1}
    \end{pmatrix}
    &=
    \begin{pmatrix}
        X^{n} \\
        V^{n}
    \end{pmatrix}
    + \frac{2\tau}{3}
    \begin{pmatrix}
        V^{n+\frac{1}{4}} \\
        \mb{f}[X^{n+\frac{1}{4}}]
    \end{pmatrix}
    + \frac{\tau}{3}
    \begin{pmatrix}
        V^{n+1} \\
        \mb{f}[X^{n+1}]
    \end{pmatrix}.
\label{E:STAGE2}
\end{align}
Equation~\eqref{E:STAGE1} is simply the backward Euler method with
steplength $\tau/4$. Therefore it can be solved as in
\eqref{E:UPDATE:LA:A}--\eqref{E:AMIN:LA} above with $\alpha$, $\tau$, $X^{n+1}$
and $V^{n+1}$
replaced by $1$, $\tau/4$, $X^{n+\frac{1}{4}}$ and $V^{n+\frac{1}{4}}$, respectively.
Next we choose a suitable linear combination of \eqref{E:STAGE1} and
\eqref{E:STAGE2} to eliminate the term $\mb{f}[X^{n+\frac{1}{4}}]$. We find
\begin{equation}
    3
    \begin{pmatrix}
        X^{n+1} \\
        V^{n+1}
    \end{pmatrix}
    -8
    \begin{pmatrix}
        X^{n+\frac{1}{4}} \\
        V^{n+\frac{1}{4}}
    \end{pmatrix}
    +5
    \begin{pmatrix}
        X^{n} \\
        V^{n}
    \end{pmatrix}
    = \tau
    \begin{pmatrix}
        V^{n+1} \\
        \mb{f}[X^{n+1}]
    \end{pmatrix}.
\label{E:RK:REDUCED}
\end{equation}
The scheme coefficients \eqref{E:BUTCHER} were chosen to minimize the
magnitude of the coefficient $8$ in \eqref{E:RK:REDUCED}. Since
$X^{n+\frac{1}{4}}$ and $V^{n+\frac{1}{4}}$ are already known from
solving \eqref{E:STAGE1}, this equation is structurally identical to
the multi-step scheme \eqref{E:BDF2:DEF}. As before, the unknown $X:=X^{n+1}$
can be found independently from $V^{n+1}$. We have
\begin{gather*}
    \frac{24}{\tau^2} \left[
        \Big( X - (X^{n+\frac{1}{4}}
            + {\TST\frac{3\tau}{4} V^{n+\frac{1}{4}}} ) \Big)
        -\frac{5}{8} \Big( X - \big( X^n
        + ( {\TST\frac{\tau}{3}V^n + \frac{2\tau}{3}V^{n+\frac{1}{4}}} )
            \big) \Big) \right]
    = \mb{f}[X],
\\
    \vphantom{\bigg|}
    V^{n+1}
        = 8 ( {\TST\frac{3}{4}V^{n+\frac{1}{4}}} )
        - 5 ( {\TST\frac{1}{3}V^n + \frac{2}{3}V^{n+\frac{1}{4}}} )
        + {\TST\frac{\tau}{3}} \mb{f}[X^{n+1}].
\end{gather*}
The first of these is solved by a minimization problem similar to
(\ref{E:AMIN:LA}) and (\ref{E:F2}) above. The left-hand side is then
substituted for $\mb{f}[X^{n+1}]$ in the second equation.  We leave
the details of spatial discretization and energy dissipation (through
a choice of new velocities that dissipate the maximum amount of energy
without changing the convected distribution of mass) to the reader as
they are similar to the BDF2 case.  The generalization to higher
order multi-step and diagonally implicit Runge-Kutta methods is also
straightforward, although we have not tested them.
\end{remark}


\subsection{Implementation details}\label{SS:ID}

The schemes of Definitions~\ref{D:SCHEME1} and \ref{D:SCHEME2} are easy
to implement once a solver for convex optimization is available. We
experimented with our own implementation of a trust-region method
described below, and with the convex optimization solvers cvxopt
\cite{cvxopt} and Knitro \cite{Knitro}. In order to illustrate the key
issues, let us consider the problem of minimizing the functional
\begin{equation}
    F(\BZ) := \sum_{i=1}^N \frac{3}{4\tau^2} m |z_i-\hat{x}^n_i|^2
        + \sum_{i=1}^{N-1} U\bigg(\frac{m}{z_{i+1}-z_i}\bigg)
            (z_{i+1}-z_i)
  \label{E:F:AGAIN}
\end{equation}
among all $\BZ\in\R^N$ such that $z_{i+1}>z_i$ for all $i$. The first
term in \eqref{E:F:AGAIN} penalizes the difference between the
positions $z_i$ and the intermediate positions $\hat{x}^n_i$, and since
$1/\tau^2$ is very large for small timesteps, it is natural to choose
the intermediate positions as the starting vector for the minimization
process. To have a bounded internal energy, we spread the $\hat{x}^n_i$
slightly apart so that the minimal distance between neighboring
positions of the starting vector is bigger than some number $dmin$. We
define
\begin{equation}
  z^\e0_i := \hat{x}^n_i + \TST\frac{1}{2}\big( l_i + r_i \big),
\end{equation}
where the negative displacements $l_i$ are built up by
\medskip
\begin{tabbing}
\hspace*{2em} \= \hspace*{2em} \= \hspace*{2em} \= \hspace*{2em} \= \kill
    \> \textbf{for} $i=2$ \textbf{to} $N$
\\
    \>\> \textbf{for} $j=i-1$ \textbf{downto} 1
\\
    \>\>\> $d=\hat{x}^n_{j+1} + l_{j+1} - (\hat{x}^n_j + l_j) - dmin$
\\
    \>\>\> \textbf{if} $d<0$
\\
    \>\>\>\> $l_j = l_j + d$
        \quad (decrease $l_j$)
\\
    \>\>\> \textbf{else}
\\
    \>\>\>\> \textbf{break}
        \quad (continue with outer loop)
\end{tabbing}
\medskip
and a similar algorithm is used to obtain the $r_i$. Note that this
procedure preserves the ordering. In practice, only a few points
$\hat{x}^n_i$ need to be adjusted and most of the resulting positions
$z_i^\e0$ are equal to $\hat{x}^n_i$. We use $dmin := \frac{1}{2}
\min_i(x_{i+1}^n-x_i^n)$.

From this initial guess, we proceed with a trust region Newton method
to minimize the functional $F(\BZ)$. Our implementation is based on
Chapter 4 of \cite{NocedalWright}. For a given $\BZ\in\R^N$, the
gradient $g=\nabla_{\BZ} F(\BZ)$ and Hessian $H=\nabla^2_{\BZ} F(\BZ)$
are trivial to compute for the first term in (\ref{E:F:AGAIN}), and are
easily assembled element by element (one interval $[z_i,z_{i+1})$ at a
time) for the second term. The Hessian is positive definite (as it is
strictly diagonally dominant) and tridiagonal.  At each step of the
minimization algorithm, we seek the minimizer of the quadratic
approximation
$$
    Q(\BP) := F(\BZ) + g^T\BP + \TST\frac{1}{2}\BP^T H\BP
        \;\approx\; F(\BZ+\BP)
$$
over all $\BP\in\R^d$ satisfying the trust region constraint
$\|D^{-1}\BP\| \LS \Delta$. Here $D$ is the diagonal scaling matrix
with entries $D_{ii} := \frac{1}{3} \min\bigl\{z_i - z_{i-1}, z_{i+1} -
z_i\bigr\}$ that prevents the positions $z_i$ from crossing. The trust
region radius $\Delta$ is not allowed to exceed one. By defining
$\hat{\BP} := D^{-1}\BP$, we map the elliptical trust region to a
sphere. We then consider the following modified problem: Minimize the
functional
\begin{equation}
    \hat{Q}(\hat{\BP}) := F(\BZ) + \hat{g}^T \hat{\BP}
        + \TST\frac{1}{2} \hat{\BP}^T \hat{H}\hat{\BP},
\label{E:QP}
\end{equation}
over all $\hat{\BP}\in\R^N$ with $\|\hat{\BP}\| \LS \Delta$. Here
$\hat{g} := Dg$ and $\hat{H} := DHD$. Notice that this minimization
problem reduces to the previous one if $D=\ID$.

The main theorem governing the design of trust region methods is

\begin{theorem}[Trust region optimality criterion]
The vector $\BP_*\in\R^N$ minimizes the quadratic functional
\eqref{E:QP} over all $\hat{\BP}\in\R^N$ with $\|\hat{\BP}\|\LS\Delta$
if and only if $\BP_*$ is feasible and there exists a scalar
$\lambda_*\GS 0$ with the following properties:
\begin{align}
    (\hat{H} + \lambda_*\,\ID)\BP_* &= -\hat{g},
\label{E:TR1}\\
    \lambda_* (\Delta - \|\BP_*\|) &= 0,
\label{E:TR2}\\
    & \hspace*{-7em}
        (\hat{H} + \lambda_*\,\ID) \text{ is positive semidefinite.}
\label{E:TR3}
\end{align}
\end{theorem}

In our case, the matrix $\hat{H}$ is positive definite, so
(\ref{E:TR3}) is satisfied automatically. Moreover, equation
(\ref{E:TR1}) has a solution for any $\lambda_*\GS 0$, which eliminates
the need to search for the most negative eigenvalue of $\hat{H}$ and
also rules out the ``hard case'' of Mor\'{e} and Sorensen; see
\cite{NocedalWright}. We use a variant of the Cholesky
factorization algorithm to find the Lagrange multiplier $\lambda_*$. In
pseudo-code it is given as
\medskip
\begin{tabbing}
\hspace*{2em} \= \hspace*{2em} \= \hspace*{2em} \= \hspace*{2em} \= \kill
    \> Given $\lambda^\e0$, $\Delta>0$
\\
    \> \textbf{for} $k=0,1,2,\dots$
\\
    \>\> Factorize $\hat{H}+\lambda^\e{k}\ID = LL^T$
\\
    \>\> Solve $LL^Tp_k = -\hat{g}$ and $Lq_k = p_k$
\\
    \>\> Define $\jd\lambda^\e{k+1} = \lambda^\e{k}
        + \left(\frac{\|p_k\|}{\|q_k\|}\right)
          \left(\frac{\|p_k\|-\Delta}{\Delta}\right)$
\\
    \>\> \textbf{if} $\lambda^\e{k}=0$ and $\lambda^\e{k+1}\le 0$
\\
    \>\>\> \textbf{return}
\\
    \>\> \textbf{if} $|(\|p_k\|-\Delta)/\Delta|<10^{-12}$
\\
    \>\>\> compute $p_{k+1}$
\\
    \>\>\> \textbf{return}
\\
    \>\> \textbf{if} $\lambda^\e{k+1}<0$
\\
    \>\>\> set $\lambda^\e{k+1}=0$
\\
    \> \textbf{end for}
\end{tabbing}
\medskip
This algorithm is equivalent to Newton's method for finding the zeros
of the function $\Delta^{-1}-\|p(\lambda)\|^{-1}$, which usually
converges to machine precision in 4--6 iterations. We vary the size
$\Delta$ of the trust region in the standard way \cite{NocedalWright}
by comparing the actual reduction of the function to the predicted
reduction of the quadratic model. For any given $\Delta$, we use the
Lagrange multiplier $\lambda_*$ of the previous trust region step as
starting value $\lambda^\e0$ for the iteration above. Since we use the
exact Hessian, the method converges quadratically (usually requiring
4-20 trust region steps). Moreover, as $\hat{H}$ is tridiagonal, the
Cholesky factorization is of complexity $O(N)$ only, and so the convex
optimization procedure can be solved very efficiently.

The minimization problem for the functional \eqref{E:F2} is almost
identical to the one we just described. The first term in
(\ref{E:F:AGAIN}) must be replaced by
$$
    \sum_{i,j=0}^N \Bigg(
        \frac{3}{2\tau^2}
            (z_i-x'_i) A_{ij} (z_j-x'_j)
        -\frac{3}{8\tau^2}
            (z_i-x''_i) A_{ij} (z_j-x''_j)
    \Bigg).
$$
The Hessian matrix for this part of the objective function is
$\frac{9}{4\tau^2}A$, which is again positive definite and tridiagonal,
and so the problem can be solved efficiently.


\section{Numerical Experiments}\label{S:NE}

In this section, we report on numerical experiments we performed with
the variational particle methods VPS1 and VPS2 introduced in
Definitions~\ref{D:SCHEME1} and~\ref{D:SCHEME2}. We studied both the
porous medium and heat equations, as well as the isentropic and
isothermal Euler equations, for various choices of parameters and
initial data.


\subsection{Porous Medium and Heat Equations}\label{SS:PMHE}

For both the porous medium equation \eqref{E:POR} and the heat equation
\eqref{E:HE}, solutions converge to self-similar profiles as
$t\rightarrow\infty$. These profiles can be computed explicitly, which
makes it possible to measure the error of the numerical approximation
exactly.


\subsubsection{Porous Medium Equation}

In the case of the porous medium equation the self-similar solution
are called Barenblatt profiles and given by the formula
\begin{equation}
    \RHO_*(t,x) = t^{-\alpha} \Big( C^2-k t^{-2\beta}|x|^2
        \Big)_+^{1/(\gamma-1)},
\label{E:BARENBLATT}
\end{equation}
where $(s)_+ := \max\{s,0\}$, and
$$
    \alpha = \frac{d}{d(\gamma-1)+2},
    \quad
    \beta = \frac{\alpha}{d},
    \quad
    k = \frac{\beta(\gamma-1)}{2\gamma}.
$$
see \cite{Vazquez}. The constant $C$ is chosen in such a way that $M(t)
:= \int_{\R^d} \RHO_*(t,x) \,dx$ equals the total mass. Notice that
$M(t)$ is in fact independent of $t>0$, and for simplicity we assume
that $M(t)=1$. Then $\RHO_*(t,\cdot)$ converges to a Dirac measure as
$t\rightarrow 0$. In the one-dimensional case under consideration, the
constants simplify to
$$
    \alpha = \frac{1}{\gamma+1},
    \quad
    \beta = \alpha,
    \quad
    k = \frac{\gamma-1}{2\gamma(\gamma+1)},
$$
and
$$
    C = \Bigg( \bigg( \frac{\gamma-1}{2\pi\gamma(\gamma+1)}
            \bigg)^{\frac{1}{2}} \;
        \frac{\Gamma\big(\frac{3}{2}+\frac{1}{\gamma-1}\big)}
            {\Gamma\big(\frac{\gamma}{\gamma-1}\big)}
        \Bigg)^{\frac{\gamma-1}{\gamma+1}},
$$
where $\Gamma(z) := \int_0^\infty t^{z-1} e^{-t} \,dt$ for all $z\GS
0$. To obtain an approximate solution of the porous medium equation
\eqref{E:POR} we used the variational particle schemes VPS1 and VPS2 of
the previous section with internal energy $U(\RHO) =
\RHO^\gamma/(\gamma-1)$.

\begin{figure}
\includegraphics[scale=0.75]{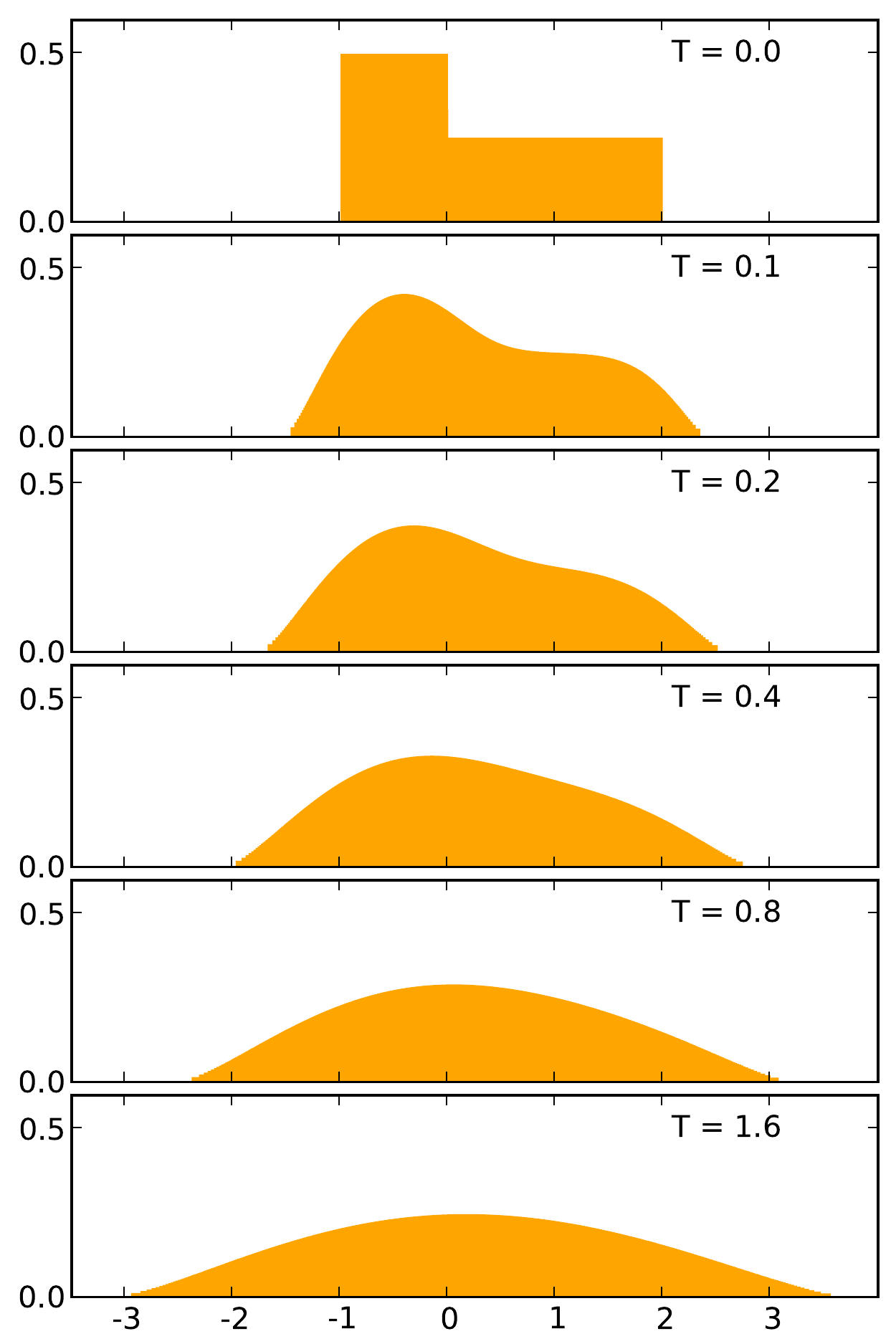}
\caption{Porous medium equation with $\gamma=5/3$.}
\label{F:PMTIME}
\end{figure}

Figure~\ref{F:PMTIME}
shows the solution of \eqref{E:POR} with $\gamma = 5/3$ at different
times. We started out with asymmetric and discontinuous initial data
$$
    \bar{\RHO}(x) := \begin{cases}
        0.5 & \text{if $x\in(-1,0)$,}
\\
        0.25 & \text{if $x\in(0,2)$,}
\\
        0 & \text{otherwise.}
    \end{cases}
$$
We used the VPS1 scheme with $N=1000$ points and $\tau=0.0016$ to
generate these plots.  The running time was 8 seconds on a 2.4 GHZ
laptop.  Note that the solution at time $T=1.6$ is close to a (shifted)
Barenblatt profile.

\begin{figure}
\includegraphics[scale=0.75]{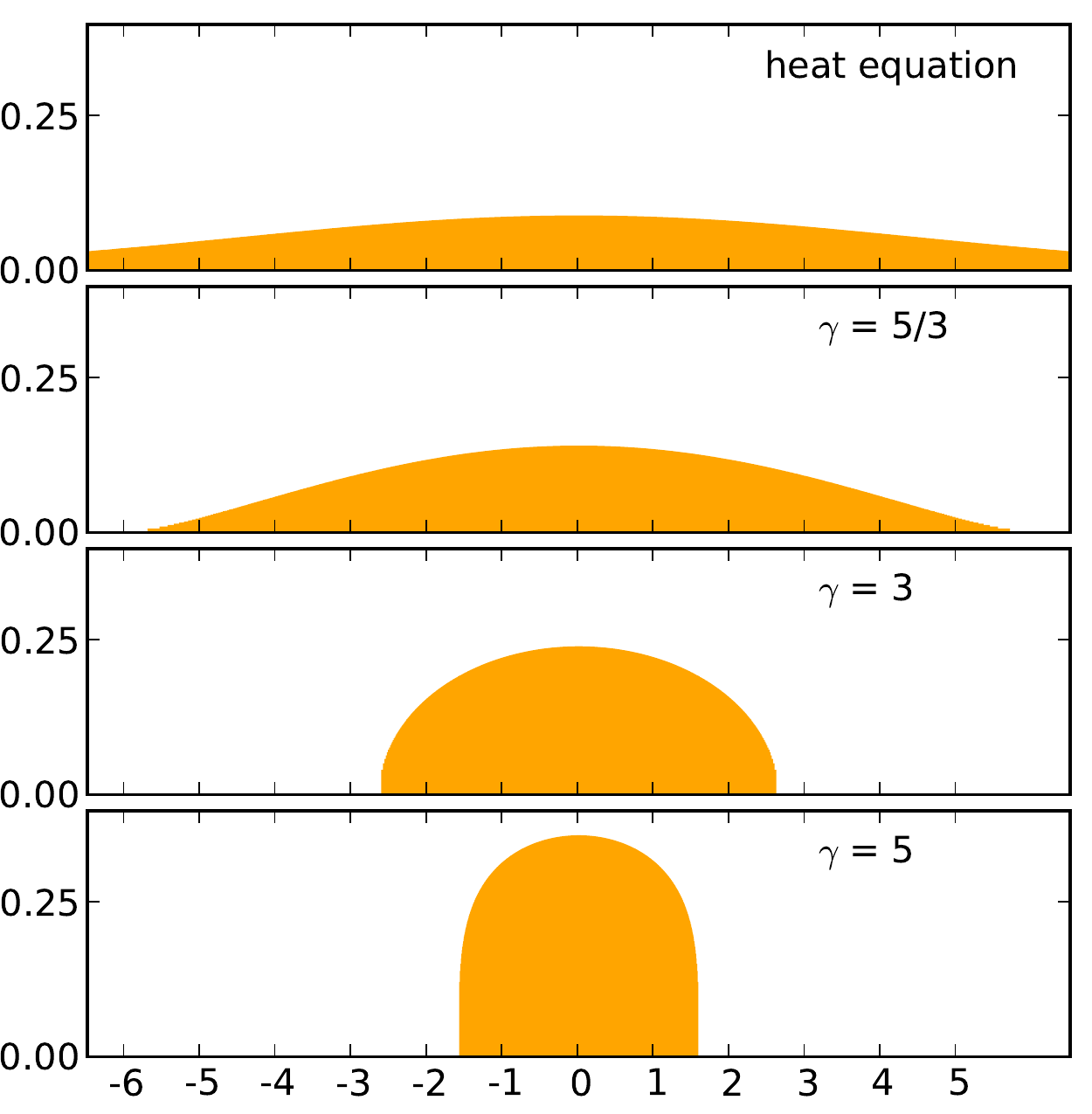}
\caption{Porous medium/heat equation at time $T=10.0$.}
\label{F:PMERROR}
\end{figure}

To estimate the accuracy of the VPS1 scheme, we used initial data
\begin{equation}
    \bar{\RHO}(x) := \begin{cases}
        50 & \text{if $x\in(-0.01,0.01)$,}
\\
        0 & \text{otherwise,}
    \end{cases}
\label{E:DIRAC}
\end{equation}
which approximates a Dirac mass, and computed the solution at time
$T=10$ for both the porous medium equation \eqref{E:POR} with different
values of $\gamma>1$, and for the heat equation, which corresponds to
the limiting case $\gamma=1$. The particle mass was chosen as $m=0.001$
and the timestep $\tau=0.01$. The results are shown in
Figure~\ref{F:PMERROR}.
%
%
Note that the solution is not continuously differentiable for
$\gamma\GS 2$: The contact angle is strictly positive for $\gamma=2$,
is vertical for $\gamma>2$, and vanishes for $\gamma<2$.
For the heat
equation, the self-similar solution is given by a Gauss function
\begin{equation}
    \RHO_*(t,x) = \frac{1}{\sqrt{4\pi t}}
        \exp\bigg(-\frac{|x|^2}{4t}\bigg).
\label{E:GAUSSIAN}
\end{equation}

The following table shows the $\L^\infty$-error of the
approximate solution:
$$
\begin{tabular}{|l||c|}
    \hline
    \text{Problem}
        & \text{$\L^\infty$-error}
\\
    \hline
    \hline
    \text{heat eq.}
        & 9.45e-5 
\\
    \hline
    $\gamma=5/3$
        & 1.53e-4 
\\
    \hline
    $\gamma=3$
        & 8.63e-4 
\\
    \hline
    $\gamma=5$
        & 2.62e-3 
\\
    \hline
\end{tabular}
$$
This error is computed, with $n:=T/\tau$, according to the
formula
\begin{equation}
    \L^\infty\text{-error}\; :=\;
 \max_{1\LS i< N} \bigg| \frac{m}{x^n_{i+1} - x^n_i}
        -\RHO_*\Big(10,\jt\frac{1}{2}(x^n_i + x^n_{i+1})\Big) \bigg|.
\label{E:EINF}
\end{equation}
Since the Barenblatt profile has steep gradients for large exponents
$\gamma$, the infinity norm of the error grows as we increase $\gamma$.
On the other hand, the support of the solution of \eqref{E:POR}
decreases as we increase $\gamma$. For the heat equation, which has an
infinite speed of propagation, the support is unbounded.  The two
effects largely cancel each other out when we estimate the $\L^1$-error
as follows:
$$
\begin{tabular}{|l||c|}
    \hline
    \text{Problem}
        & \text{$\L^1$-error}
\\
    \hline
    \hline
    \text{heat eq.}
        & 1.11e-3 
\\
    \hline
    $\gamma=5/3$
        & 1.08e-3  
\\
    \hline
    $\gamma=3$
        & 9.62e-4  
\\
    \hline
    $\gamma=5$
        & 7.26e-4  
\\
    \hline
\end{tabular}
$$
Here the error is computed, with $n:=T/\tau$, according to the
formula
\begin{equation}
    \L^1\text{-error}\; := \;
    \sum_{i=1}^{N-1} \bigg| \frac{m}{x^n_{i+1} - x^n_i}
        -\RHO_*\Big(10,\jt\frac{1}{2}(x^n_i + x^n_{i+1})\Big)
            \bigg| (x^n_{i+1}-x^n_i).
\label{E:EONE}
\end{equation}

To compute the convergence rate of the VPS1 method, we computed the
solution at time $T=10$ of the porous medium equation with initial data
\eqref{E:DIRAC} and $\gamma=5/3$, for several choices of particle mass
$m$ and timestep $\tau$. We have
$$
\begin{tabular}{|c|c||c|c||c|c||c|c|}
    \hline
    &
    & \text{Error at}
    &
    &
    &
    &
    &
\\[-0.5ex]
    \text{$m$}
    & \text{$\tau$}
    & \text{$x=0$}
    & \text{Rate}
    & \text{$\L^\infty$-error}
    & \text{Rate}
    & \text{$\L^1$-error}
    & \text{Rate}
\\
    \hline
    \hline
    0.01   & 0.1    &  1.15e-3 &      & 1.15e-3 &      & 6.79e-3 &
\\
    0.004  & 0.04   &  5.21e-4 & 0.86 & 5.21e-4 & 0.86 & 3.35e-3 & 0.77
\\
    0.001  & 0.01   &  1.53e-4 & 0.88 & 1.53e-5 & 0.88 & 1.08e-3 & 0.82
\\
    0.0004 & 0.004  &  6.73e-5 & 0.90 & 6.73e-5 & 0.90 & 4.98e-4 & 0.84
\\
    0.0001 & 0.001  &  1.90e-5 & 0.91 & 3.34e-5 & 0.50 & 1.47e-4 & 0.88
\\
    \hline
\end{tabular}
$$
Thus, the VPS1 method performs slightly worse than a first order method
for these initial conditions. The last data point for the
$\L^\infty$-norm is anomalous because the location of the largest error
jumps from the center $x=0$ to the boundary of the support of the
solution, which is not a smooth transition.

In this experiment, the initial data was rather singular. To estimate
the effect of this irregularity, we also performed a computation
starting off from the Barenblatt profile at time $t=1$, which is much
more regular, and continued to $t=2$. More specifically, we partitioned
the support $[-C/\sqrt{k},C/\sqrt{k}]$ of $\RHO_*(1,\cdot)$ into
intervals $I_i=[\tilde{x}_{i-1}, \tilde{x}_i]$ such that $\int_{I_i}
\RHO_*(1,x)\,dx=1/N$ for all $1\LS i\LS N$, and define $x_i$ to be the
center of mass of $\RHO_*(1,\cdot)$ over $I_i$.  Again we used
$\gamma=5/3$. We have
$$
\begin{tabular}{|c|c||c|c||c|c||c|c|}
    \hline
    &
    & \text{Error at}
    &
    &
    &
    &
    &
\\[-0.5ex]
    \text{$m$}
    & \text{$\tau$}
    & \text{$x=0$}
    & \text{Rate}
    & \text{$\L^\infty$-error}
    & \text{Rate}
    & \text{$\L^1$-error}
    & \text{Rate}
\\
    \hline
    \hline
    0.01   & 0.05    &  4.32e-4  &      &  7.99e-4   &     &  1.58e-3   &
\\
    0.004  & 0.02    &  1.88e-4  &  0.91 &  5.69e-4   & 0.37 &  7.88e-4   & 0.76
\\
    0.001  & 0.005   &  5.00e-5  &  0.96 &  2.89e-4   & 0.49 &  2.38e-4   & 0.86
\\
    0.0004 & 0.002   &  2.04e-5  &  0.98 &  1.76e-4   & 0.54 &  1.03e-4   & 0.91
\\
    0.0001 & 0.0005  &  5.17e-6  &  0.99 &  8.03e-5  & 0.57 &  2.80e-5  & 0.94
\\
    \hline
\end{tabular}
$$
In the interior of the support of the solution, the method converges at
first order, but near the boundary the order deteriorates. The
$\L^1$-error involves a mixture of the two convergence rates, and is
therefore slightly worse than first order.

We also tested the convergence rates of the VPS2 method.  For initial
data of the form (\ref{E:DIRAC}) with $\gamma=5/3$, we chose $N+1$
knots $x_i$ uniformly spaced over the interval $[-0.01,0.01]$ and set
$m_i = 1/N$ for all $1\LS i\LS N$.  The errors are smaller than they
were for VPS1, but the method still converges at first order only due
to the discontinuity in the initial density. We have
$$
\begin{tabular}{|c|c||c|c||c|c||c|c|}
    \hline
    &
    & \text{Error at}
    &
    &
    &
    &
    &
\\[-0.5ex]
    \text{$m$}
    & \text{$\tau$}
    & \text{$x=0$}
    & \text{Rate}
    & \text{$\L^\infty$-error}
    & \text{Rate}
    & \text{$\L^1$-error}
    & \text{Rate}
\\
    \hline
    \hline
    0.01   & 0.1    & 1.13e-4  &      & 1.34e-3  &      & 2.34e-3 &
\\
    0.004  & 0.04   & 4.37e-5  & 1.04 & 8.46e-4  & 0.50 & 1.02e-3 & 0.90
\\
    0.001  & 0.01   & 1.06e-5  & 1.02 & 3.97e-4  & 0.55 & 2.84e-4 & 0.93
\\
    0.0004 & 0.004  & 4.18e-6  & 1.01 & 2.36e-4  & 0.57 & 1.19e-4 & 0.95
\\
    0.0001 & 0.001  & 1.03e-6  & 1.01 & 1.05e-4  & 0.58 & 3.12e-5 & 0.97
\\
    \hline
\end{tabular}
$$
To estimate the influence of the singularity of the initial data, we
tried the VPS2 method starting again off from the Barenblatt profile at
time $t=1$ and continuing to $t=2$. We partitioned the support
$[-C/\sqrt{k},C/\sqrt{k}] =: [a,b]$ into equal subintervals with knots
$x_i = a + i(b-a)/N$ for all $1\LS i\LS N$ and set
\begin{equation}
    m_i := \int_{x_{i-1}}^{x_i} \RHO_*(1,x)\,dx
    \quad\text{for all $1\LS i\LS N$.}
\label{E:mi}
\end{equation}
We used $\gamma=5/3$. When comparing the result to $\RHO_*(2,\cdot)$,
we have
$$
\begin{tabular}{|c|c||c|c||c|c||c|c|}
    \hline
    &
    & \text{Error at}
    &
    &
    &
    &
    &
\\[-0.5ex]
    \text{$N$}
    & \text{$\tau$}
    & \text{$x=0$}
    & \text{Rate}
    & \text{$\L^\infty$-error}
    & \text{Rate}
    & \text{$\L^1$-error}
    & \text{Rate}
\\
    \hline
    \hline
    100   & 0.1    & 1.64e-5  &      &  3.44e-4  &       &  1.02e-4  &
\\
    250   & 0.04   & 2.64e-6  & 1.99 &  9.82e-5  & 1.37  &  1.71e-5  & 1.95
\\
    1000  & 0.01   & 1.66e-7  & 2.00 &  1.32e-5  & 1.45  &  1.10e-6  & 1.98
\\
    2500  & 0.004  & 2.66e-8  & 2.00 &  3.40e-6  & 1.48  &  1.78e-7  & 1.99
\\
    10000 & 0.001  & 1.66e-9  & 2.00 &  4.30e-7  & 1.49  &  1.13e-8  & 1.99
\\
    \hline
\end{tabular}
$$
As expected, we obtain second order convergence throughout most of the
interval when we start with continuously differentiable initial data.
A few points next to the boundary of the support, however, seem to
converge at a lower rate of $3/2$. This raises the question of what
happens if the initial data is continuous but not $\C^1$.  The
Barenblatt profile with $\gamma=3$ has a square-root singularity at the
edges of its support. Repeating the above procedure for this case, we
obtain
$$
\begin{tabular}{|c|c||c|c||c|c||c|c|}
    \hline
    &
    & \text{Error at}
    &
    &
    &
    &
    &
\\[-0.5ex]
    \text{$N$}
    & \text{$\tau$}
    & \text{$x=0$}
    & \text{Rate}
    & \text{$\L^\infty$-error}
    & \text{Rate}
    & \text{$\L^1$-error}
    & \text{Rate}
\\
    \hline
    \hline
    100   & 0.1    & 2.36e-5  &      & 5.49e-4  &      &  1.22e-4  &
\\
    250   & 0.04   & 3.78e-6  & 2.00 & 3.43e-4  & 0.51 &  3.16e-5  & 1.48
\\
    1000  & 0.01   & 2.37e-7  & 2.00 & 1.67e-4  & 0.52 &  4.02e-6  & 1.49
\\
    2500  & 0.004  & 3.80e-8  & 2.00 & 1.04e-4  & 0.51 &  1.02e-6  & 1.50
\\
    10000 & 0.001  & 2.37e-9  & 2.00 & 5.17e-5  & 0.51 &  1.28e-7  & 1.50
\\
    \hline
\end{tabular}
$$
It appears that the infinite slope of the exact solution near the
boundary of the support creates large errors there that dominate the
overall $\L^1$-error of the method. The situation can be improved if we
redistribute the mass to better resolve the solution near these
endpoints. Specifically, we choose the knot positions as
$$
  x_i := f\bigg(-1+\frac{2i}{N}\bigg) \frac{C}{\sqrt{k}}
  \quad\text{for all $0\LS i\LS N$,}
$$
with weight function
$$
    f(x) := \frac{\DST\int_0^x \sqrt{1-y^2}\,dy}
        {\DST\int_0^1\sqrt{1-y^2}\,dy}
    \quad\text{for all $x\in[-1,1]$,}
$$
and we assign to each interval $[x_{i-1},x_i)$ the mass $m_i$ defined
in (\ref{E:mi}). We have
$$
\begin{tabular}{|c|c||c|c||c|c||c|c|}
    \hline
    &
    & \text{Error at}
    &
    &
    &
    &
    &
\\[-0.5ex]
    \text{$N$}
    & \text{$\tau$}
    & \text{$x=0$}
    & \text{Rate}
    & \text{$\L^\infty$-error}
    & \text{Rate}
    & \text{$\L^1$-error}
    & \text{Rate}
\\
    \hline
    \hline
    100   & 0.1    &  2.36e-5  &     &  2.84e-3  &      &  1.20e-4  &
\\
    250   & 0.04   &  3.79e-6  & 2.0 &  1.58e-3  & 0.64 &  2.22e-5  & 1.85
\\
    1000  & 0.01   &  2.38e-7  & 2.0 &  6.07e-4  & 0.69 &  1.62e-6  & 1.89
\\
    2500  & 0.004  &  3.80e-8  & 2.0 &  3.14e-4  & 0.72 &  2.79e-7  & 1.92
\\
    10000 & 0.001  &  2.38e-9  & 2.0 &  1.14e-4  & 0.73 &  1.89e-8  & 1.94
\\
    \hline
\end{tabular}
$$
In spite of the singularity in slope near the endpoints of the support,
we achieve second order accuracy over enough of the interval to
converge at second order in the $\L^1$-norm. The $\L^\infty$-norm is
larger because the first and last intervals are smaller than the
interior intervals, so their midpoints are closer to the singularity.
The first function we tried for the mass redistribution, namely
$f(x)=\frac{3}{2}x-\frac{1}{2}x^3$, caused difficulties for the convex
optimization solver for $N=10000$.


\subsubsection{Heat Equation}\label{SSS:HE}

Both schemes perform similarly for the heat equation \eqref{E:HE} as
they did for the porous medium equation. Here we report only the VPS2
results. Using the internal energy $U(\RHO) = \RHO\log\RHO$ and initial
data \eqref{E:DIRAC}, and comparing the solution at $t=10$ to the
Gaussian $\RHO_*(10,\cdot)$ defined in \eqref{E:GAUSSIAN}, we find
$$
\begin{tabular}{|c|c||c|c||c|c||c|c|}
    \hline
    &
    & \text{Error at}
    &
    &
    &
    &
    &
\\[-0.5ex]
    \text{$m$}
    & \text{$\tau$}
    & \text{$x=0$}
    & \text{Rate}
    & \text{$\L^\infty$-error}
    & \text{Rate}
    & \text{$\L^1$-error}
    & \text{Rate}
\\
    \hline
    \hline
    0.01   & 0.1    &  1.42e-4  &      &  5.86e-4  &      &  5.12e-3  &
\\
    0.004  & 0.04   &  5.42e-5  & 1.05 &  2.73e-4  & 0.83 &  2.24e-3  & 0.90
\\
    0.001  & 0.01   &  1.27e-5  & 1.05 &  8.07e-5  & 0.88 &  6.13e-4  & 0.93
\\
    0.0004 & 0.004  &  4.78e-6  & 1.06 &  3.52e-5  & 0.91 &  2.54e-4  & 0.96
\\
    0.0001 & 0.001  &  1.02e-6  & 1.12 &  9.78e-6  & 0.92 &  6.47e-5  & 0.99
\\
    \hline
\end{tabular}
$$
As before, we obtain first order convergence for the VPS2 scheme for
discontinuous initial data. Note that we could obtain a more precise
error estimate by comparing the numerical solution to the exact
solution of \eqref{E:HE} with data \eqref{E:DIRAC}, which can be
computed explicitly. But the $\L^\infty$-distance between this exact
solution and $\RHO_*(10,\cdot)$ is around $7.4\times 10^{-8}$, which is
small in comparison to the numerical errors.

Finally, we use the VPS2 scheme to evolve the heat kernel
$\RHO_*(1,\cdot)$ to time $t=2$ via the heat equation \eqref{E:HE}. As
the support of $\RHO_*(1,\cdot)$ is infinite, instead of fixing the
knot positions first, we start by distributing the mass according to
$$
    m_i := \frac{\DST f\bigg(\frac{i}{N+1}\bigg)}
        {\DST\sum_{i=1}^N f\bigg(\frac{i}{N+1}\bigg)}
    \quad\text{for all $1\LS i\LS N$,}
$$
with weight function
$$
    f(x) := q(x)q(1-x)
    \quad\text{and}\quad
    q(x) := 10x^2+x/10.
$$
We then initialize the knot positions $x^0_i$ for all $1\LS i\LS N-1$
as
$$
  x_i^0 := -2\opn{erfc}^{-1}(2s_i)
  \quad\text{with}\quad
  s_i := \sum_{j=1}^i m_j.
$$
This construction implies that $m_i = \int_{x_{i-1}}^{x_i}
\RHO_*(1,x)\,dx$ for all $2\le i\le N-1$. Since the mass of the first
and last interval are small, any reasonable choice of $x_0$ and $x_N$
works well. We used $x_0^0 := 3x_1^0-2x_2^0$ and $x_N^0 := 3x_{N-1}^0 -
2x_{N-2}^0$. We have
$$
\begin{tabular}{|c|c||c|c||c|c||c|c|}
    \hline
    &
    & \text{Error at}
    &
    &
    &
    &
    &
\\[-0.5ex]
    \text{$N$}
    & \text{$\tau$}
    & \text{$x=0$}
    & \text{Rate}
    & \text{$\L^\infty$-error}
    & \text{Rate}
    & \text{$\L^1$-error}
    & \text{Rate}
\\
    \hline
    \hline
    100   & 0.1    &  7.02e-6   &      &  1.87e-5  &      &  1.82e-4  &
\\
    250   & 0.04   &  1.10e-6   & 2.02 &  3.04e-6  & 1.98 &  3.13e-5  & 1.92
\\
    1000  & 0.01   &  6.82e-8   & 2.01 &  1.91e-7  & 2.00 &  2.10e-6  & 1.95
\\
    2500  & 0.004  &  1.09e-8   & 2.00 &  3.05e-8  & 2.00 &  3.50e-7  & 1.96
\\
    10000 & 0.001  &  6.80e-10  & 2.00 &  1.90e-9  & 2.00 &  2.31e-8  & 1.96
\\
    \hline
\end{tabular}
$$
Thus, for the heat equation, even the $\L^\infty$-norm converges at
second order as there is no singularity in the slope of the exact
solution over the course of this simulation. The $\L^1$-norm converges
slightly slower than second order due to the fact that as we add
points, the support of the numerical solution grows.

In conclusion, we have found our variational particle scheme to be an
effective method of solving the porous medium and heat equations.
Moreover, the second order version of the method does indeed give
second order accuracy as long as the initial conditions are not too
singular.


\subsection{Isentropic Euler Equations}\label{SS:IEE}

As few smooth solutions of the isentropic Euler equations \eqref{E:IE}
are known explicitly, we test our numerical methods for
self-consistency using smooth initial data for a short enough time that
shocks do not occur.  To study their ability to capture shocks and
rarefaction waves, we test our schemes using piecewise constant initial
data, for which the Riemann problems can be solved analytically.  In
general, we find that shocks cause the second order method to perform
as a first order method, and the first order method to perform slightly
worse.  We mention in passing that these schemes also work well for the
isothermal equations, but we did not include examples for the sake of
brevity.

\subsubsection{Smooth solutions}

\begin{figure}
\includegraphics[width=.9\linewidth]{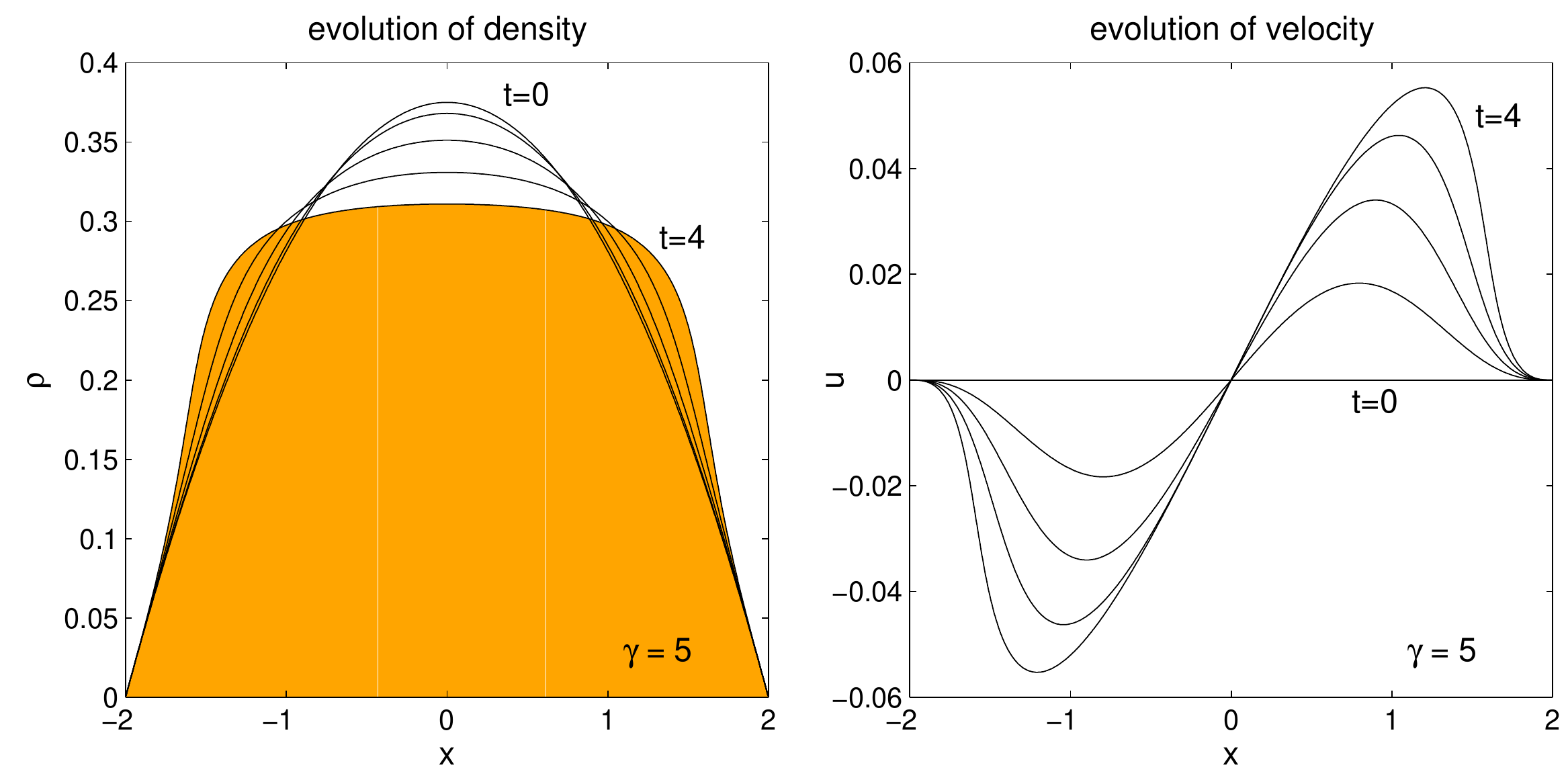} \\
\includegraphics[width=.9\linewidth]{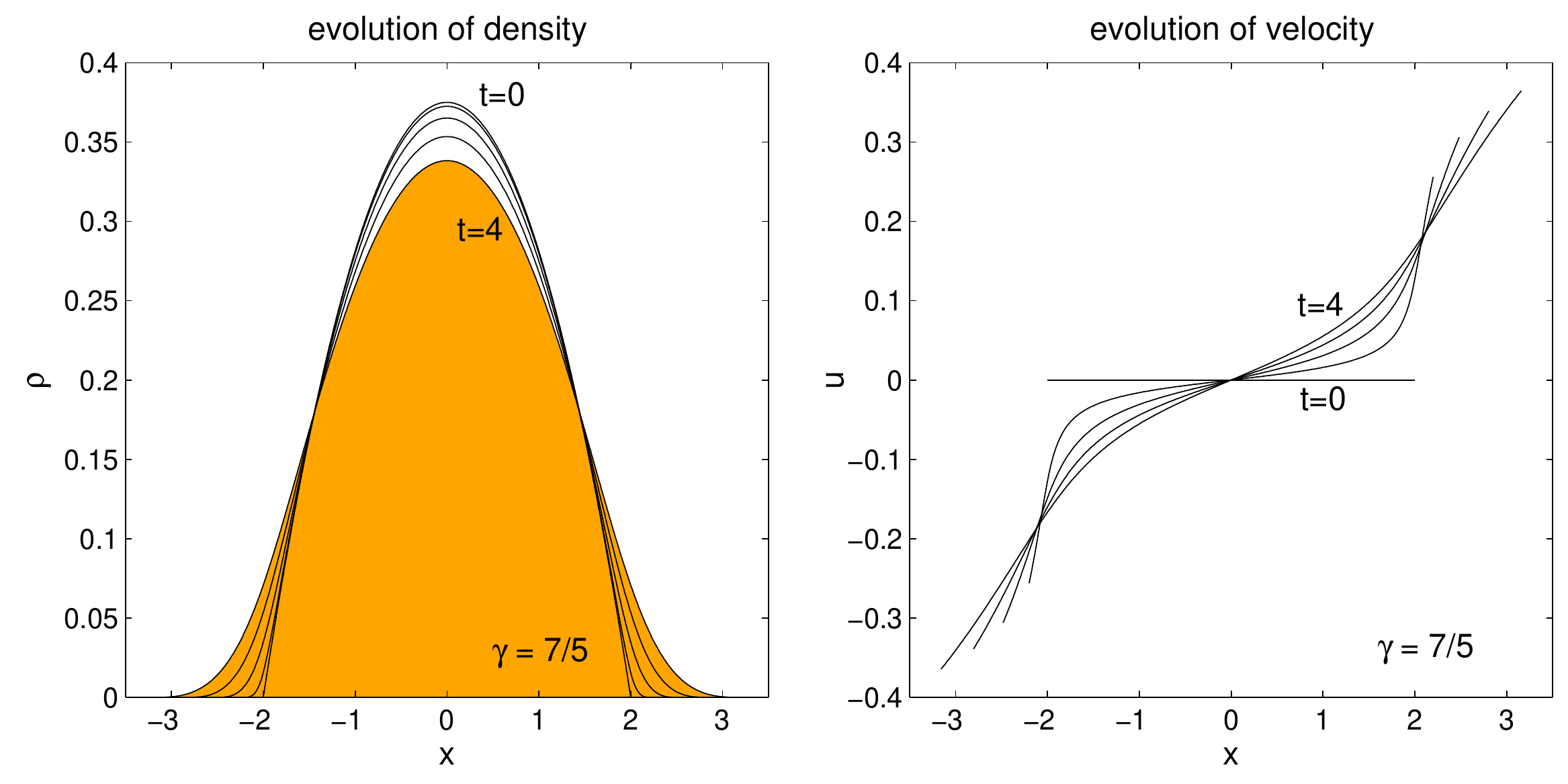}
\caption{Solution of isentropic Euler equations with a parabolic
initial distribution of mass, starting from rest. ($N=800$, VPS2)
}
\label{F:SMOOTH}
\end{figure}

We begin by testing the convergence rates of the VPS1 and VPS2 schemes
for the isentropic Euler equations with initial data
$$
    \bar\RHO(x) = \frac{3}{8} \bigg(1-\frac{x^2}{4}\bigg)_+
    \quad\text{and}\quad
    \bar{u}(x) = 0
    \quad\text{for all $x\in\R$,}
$$
for $\gamma=5$ and $\gamma=7/5$. The solutions are shown in
Figure~\ref{F:SMOOTH} at times $t=0,1,2,3,4$.
For the VPS1 scheme, we choose initial positions $x^0_i$ such that
$$
    \int_{-2}^{x^0_i} \bar{\RHO}(x) \,dx = \frac{i-1/2}{N}
    \quad\text{for all $1\LS i\LS N$.}
$$
For the VPS2 scheme with $\gamma=5$, we choose initial position $x^0_i$
uniformly spaced in the support $[-2,2]$ of the initial density, and
define masses
\begin{equation}
    m_i := \int_{x^0_{i-1}}^{x^0_i} \bar{\RHO}(x,0) \,dx
    \quad\text{for all $1\LS i\LS N$.}
\label{E:MI}
\end{equation}
For the VPS2 scheme with $\gamma=7/5$, it is necessary to distribute
more points near the ends of the support to achieve full second order
accuracy. We found that
$$
    x_i := 2f(-1+2i/N)
    \quad\text{for all $0\LS i\LS N$,}
$$
with weight function
$$
    f(x) := \frac{\DST\int_0^x \sqrt{1-y^2} \,dy}
        {\DST\int_0^1\sqrt{1-y^2}\,dy}
    \quad\text{for all $x\in[-1,1]$,}
$$
and with $m_i$ defined in (\ref{E:MI}) works well.  We used the VPS2
solution with $N=4000$ and $\tau=0.004$ as the ``exact'' solution to
compute the errors for both schemes.
$$
\begin{tabular}{|c|c||c|c||c|c|}
    \hline
    \multicolumn{2}{|c||}{} &
    \multicolumn{2}{c||}{$\gamma=5$} &
    \multicolumn{2}{c|}{$\gamma=7/5$} \\ \hline
    $N$
    & $\tau$
    & VPS1
    & VPS2
    & VPS1
    & VPS2
    \\
    \hline
    \hline
    50    & 0.4    &  1.12e-2  & 1.46e-3  &  7.75e-3  &  3.07e-3 \\
    100   & 0.2    &  8.83e-3  & 3.90e-4  &  4.68e-3  &  8.38e-4 \\
    200   & 0.1    &  4.57e-3  & 1.01e-4  &  2.51e-3  &  2.18e-4 \\
    400   & 0.05   &  2.37e-3  & 2.61e-5  &  1.24e-3  &  5.54e-5 \\
    800   & 0.025  &  1.19e-3  & 6.58e-6  &  5.84e-4  &  1.37e-5 \\
    1600  & 0.0125 &  5.96e-4  & 1.54e-6  &  3.08e-4  &  3.14e-6 \\
    \hline
\end{tabular}
\,\fbox{
\parbox[c][1.285in][b]{1.33in}{
\;\includegraphics[width=1.3in]{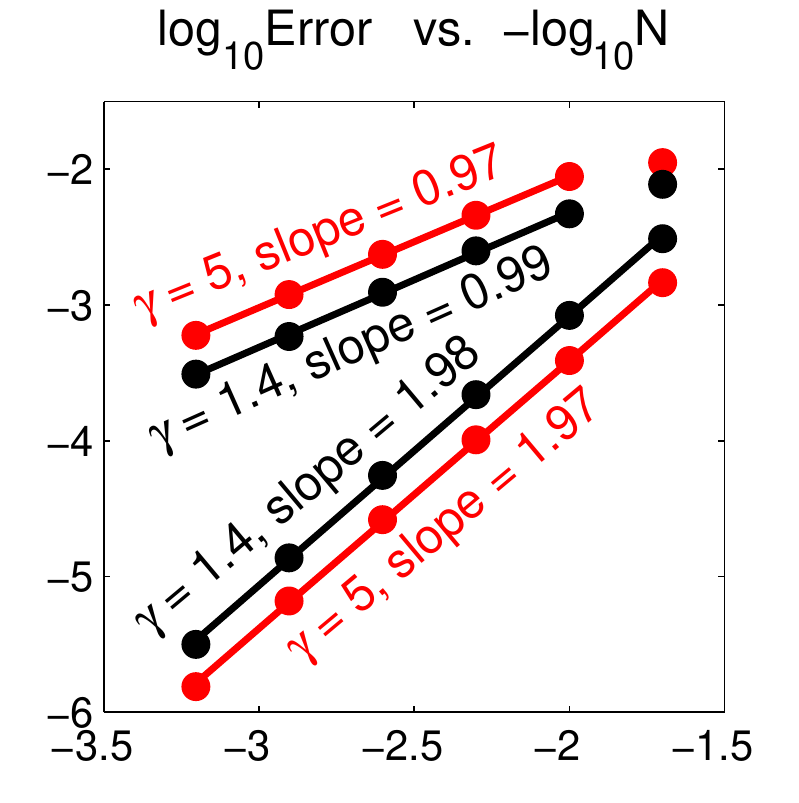}
}}
$$
The errors reported for the VPS1 scheme are estimates of the
$\L^\infty$-error of $\RHO$ alone, measured at midpoints of the
numerical solution against the ``exact'' solution, interpolated
linearly between its own midpoints. The errors reported for the VPS2
scheme involve the Wasserstein distance of the densities and the
difference in velocities. Notice that the finiteness of the total
energies only implies that
$$
    u^n \in \L^2(\R,\RHO^n)
    \quad\text{and}\quad
    u^\mathrm{exact} \in \L^2(\R^d, \RHO^\mathrm{exact}),
$$
so in general the two velocities are in different spaces. Let $\BR^n$
and $\BR^\mathrm{exact}$ denote the inverse distribution functions of
the measures $\RHO^n\LEB^1$ and $\RHO^\mathrm{exact}\LEB^1$ resp. As
explained in Section~\ref{SS:SOS}, we then have $\BR^n \#
(\IND_{[0,1]}\LEB^1) = \RHO^n\LEB^1$, which implies that
$$
    \int_\R |u^n|^2 \RHO^n \,dx
        = \int_\R |u^n|^2 \,\Big( \BR^n \# (\IND_{[0,1]}\LEB^1)
            \Big)(dx)
        = \int_{[0,1]} |u^n\circ\BR^n|^2 \,ds.
$$
That is, the pull-back $u\circ\BR^n$ of the velocity field $u^n$ under
the map $\BR^n$ is an element in $\L^2([0,1])$. The same is true for
$u^\mathrm{exact}\circ \BR^\mathrm{exact}$. We therefore define
\begin{equation}
  E_\WAS := \sqrt{
    \WAS\big(\RHO^n,\RHO^\text{exact}(\cdot,t^n)\big)^2
        + \frac{1}{2}\int_{[0,1]} \big| u^n\circ\BR^n
            - u^\text{exact}\circ\BR^\text{exact}(\cdot,t^n)
                \big|^2 \,ds },
\label{E:EW:DEF}
\end{equation}
where $t^n=4$ is the final time. Since convergence in the Wasserstein
space implies only weak* convergence in the sense of mesaures (see
\cites{AmbrosioGigliSavare, Villani}), our error estimate is weaker
that one that involves the $\L^1$-norm. The $\L^\infty$-error of the
densities of the VPS2 solutions are about three times larger than the
$E_\WAS$-errors when $\gamma=5$, and about three times smaller than the
$E_\WAS$-errors when $\gamma=7/5$.
\medskip

In the following three sections, we study the performance of our scheme
on problems involving various combinations of shocks and rarefaction
waves emanating from an initial discontinuity at the origin. In the
presence of shocks, it no longer makes sense to use the
$\L^\infty$-norm to measure the error.


\subsubsection{Shock/Shock}

\begin{figure}
\includegraphics[width=.8\linewidth]{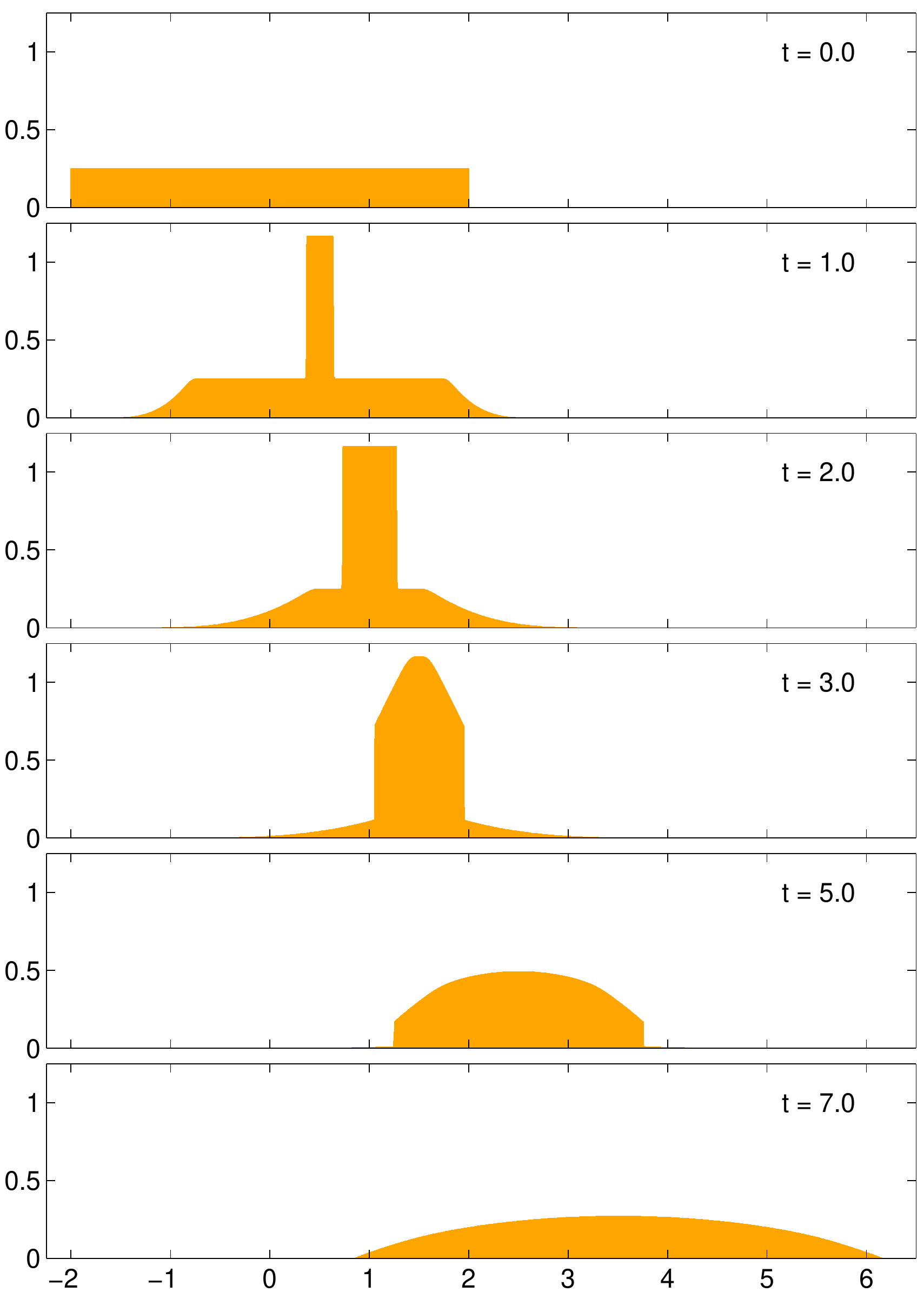}
\caption{Isentropic Euler equations with $\gamma=5/3$.}
\label{F:SSTIME}
\end{figure}

We use $\gamma=5/3$ and initial data
\begin{equation}
    (\bar{\RHO},\bar{u}) = \begin{cases}
        (0.25, 1) & \text{if $x\in(-2,0)$,}
\\
        (0.25, 0) & \text{if $x\in(0,2)$,}
\\
        (0, 0) & \text{otherwise.}
    \end{cases}
\label{E:SSDATA}
\end{equation}
Note that this amounts to solving three Riemann problems: one at
$x=0$, and two others at $x=\pm2$. Figure~\ref{F:SSTIME}
shows the evolution at different times, computed via a first order
version of the VPS2 method (with variable mass
$$
    m_i := \frac{6}{N}\int_{i-1}^i \bigg(\frac{x}{N}\bigg)
        \bigg(1-\frac{x}{N}\bigg)\,dx
    \quad\text{for all $1\LS i\LS N$,}
$$
but a first order time-stepper), with $N=1000$ and timestep
$\tau=0.01$. We will refer to this hybrid method as VPS1a. Since there
cannot be any shocks connected to the vacuum, the solution immediately
forms two rarefaction waves at the boundaries of the support that grow
in width in a self-similar fashion. The Riemann problem at $x=0$
evolves into a self-similar, double shock structure with constant
intermediate state $(\RHO_m,u_m) = (1.16641,0.5)$ with large density.
For simplicity, we only plot the density here, omitting the velocity.
As time moves on, the different waves eventually interact and the shock
strengths decrease. At time $t=7$, the solution has developed into a
continuous profile (not unlike the Barenblatt profiles of the porous
medium equation), which seems to evolve in a self-similar way while
moving to the right with a constant background speed.

\begin{figure}
\includegraphics[width=.95\linewidth]{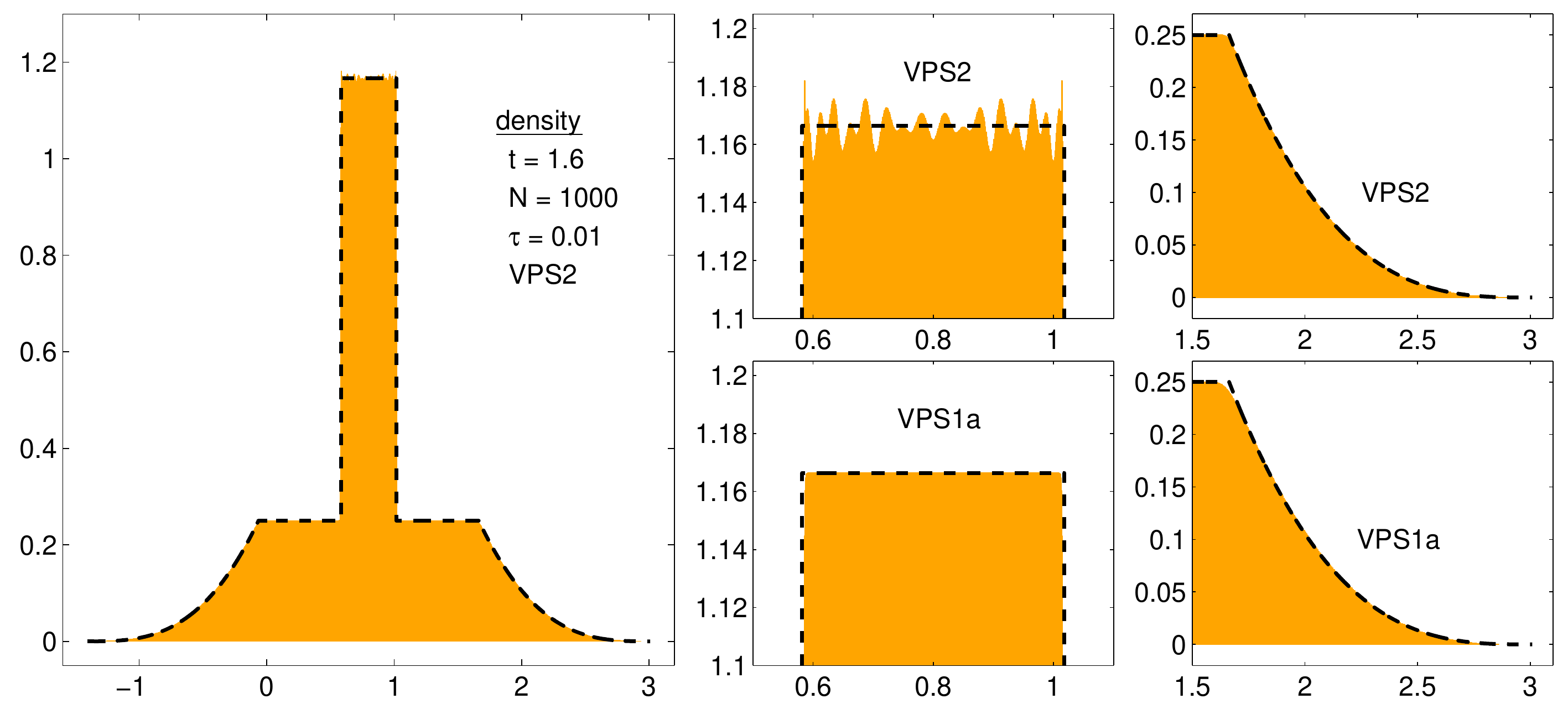}
\includegraphics[width=.95\linewidth]{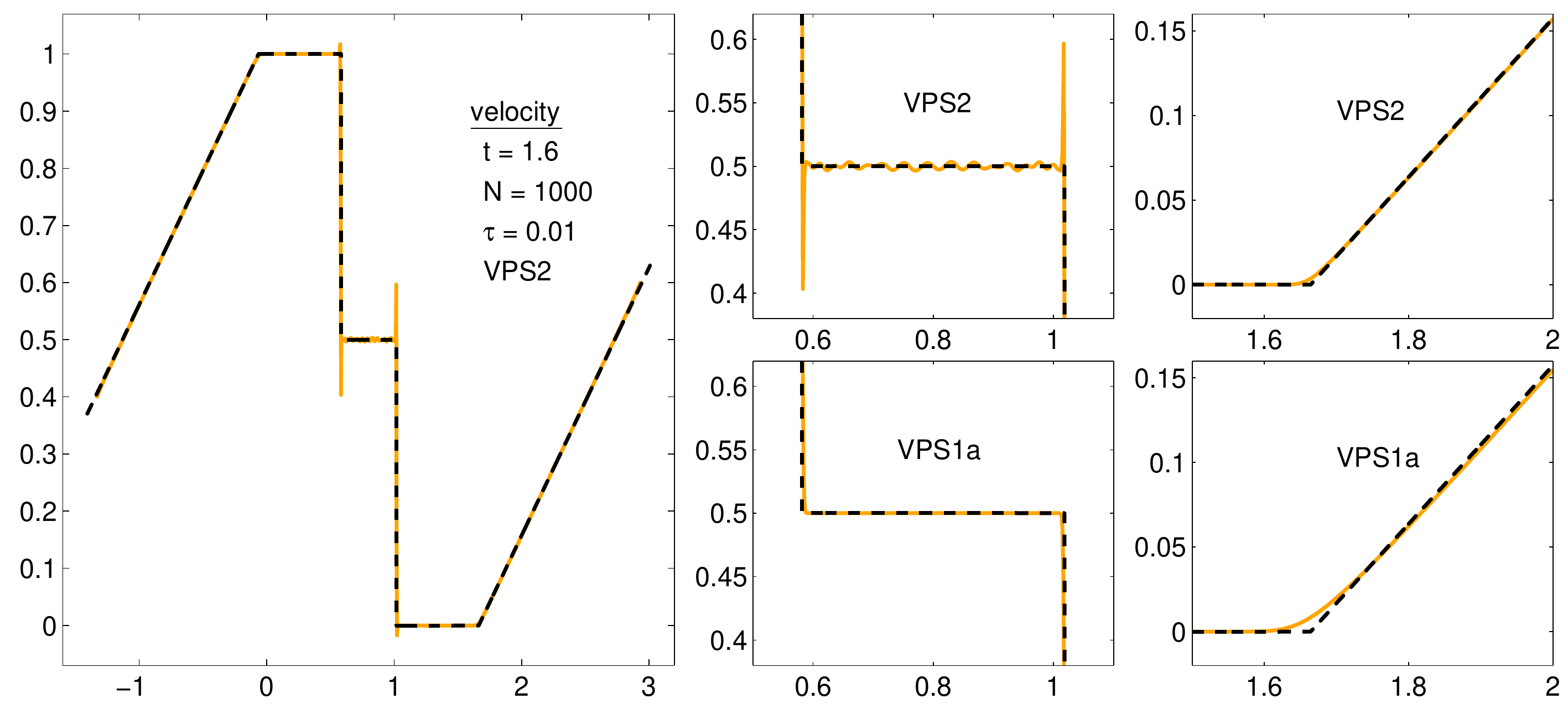}
\caption{Comparison of the VPS2 and VPS1a results with the exact solution
of the Riemann problems (dashed lines).}
\label{F:CHATTER}
\end{figure}

The VPS1 and VPS2 schemes give similar results to this hybrid method,
but the former is less accurate inside the rarefaction waves (due to
lack of resolution) while the latter suffers from small oscillations
near shocks.  These oscillations are rather interesting as they do not
grow exponentially in time nor prevent the scheme from converging as we
refine the mesh; see Figure~\ref{F:CHATTER}.
Instead, they arise because once we choose a spatial discretization
($N$ together with the mass $m_i$ of each interval), the time-stepper
is effectively solving a coupled system of masses connected by springs.
Indeed, except for the step in which mass is re-distributed when the
particles cross, we are solving Lagrange's equations of motion for the
particle trajectories $x_i(t)$, with the kinetic and potential energy
of the system defined in terms of the $x_i$ and $u_i=\dot{x}_i$ in the
obvious way (depending on the scheme).  If we were to take the
time-step $\tau$ to zero holding $N$ constant, the particles would
never cross during the transport step and the scheme would convert the
energy that is supposed to be dissipated by the shock into lattice
vibrations.  However, by simultaneously refining the mesh as we decrease
the time-step (with $\tau\sim N^{-1}$), the backward-Euler and
BDF2 methods seem to damp out these vibrations without significantly
smoothing the shocks, which is somewhat amazing.  We emphasize that
although the VPS1a scheme does not suffer from oscillations for the
timestep used in Figure~\ref{F:CHATTER}, they will eventually appear if
we take $\tau\rightarrow0$ holding $N$ constant.

Next we wish to use the exact solution of the Riemann problems to
compute the errors of the schemes and check convergence rates.
For sufficiently small $t>0$, the exact density is given by
$$
    \RHO(x,t) = \begin{cases}
        \bigg( \DST\frac{(-u_l+\RHO_l^\theta) + (x-x_l)/t}{\theta+1}
                \bigg)^{1/\theta}
            & \text{if $x\in (x_1, x_2)$,} \\
        \RHO_l & \text{if $x\in (x_2, x_3)$,} \\
        \RHO_m & \text{if $x\in (x_3, x_4)$,} \\
        \RHO_r & \text{if $x\in (x_4, x_5)$,} \\
        \bigg( \DST\frac{(u_r+\RHO_r^\theta) - (x-x_r)/t}{\theta+1}
                \bigg)^{1/\theta}
            & \text{if $x\in (x_5, x_6)$,} \\
        0 & \text{otherwise,}
    \end{cases}
$$
with constants
\begin{align*}
    x_1 &= x_l + t (u_l-\RHO_l^\theta),  &
    x_3 &= t s_l, &
    x_5 &= x_r + t (u_r-\theta\RHO_r^\theta),
\\
    x_2 &= x_l + t (u_l+\theta\RHO_l^\theta), &
    x_4 &= t s_r, &
    x_6 &= x_r + t (u_r+\RHO_r^\theta);
\end{align*}
see \cite{ChenWang}. For the case under consideration, we have
\begin{alignat}{3}
    x_l &= -2, &
    \qquad\quad \RHO_l &= 0.25, &
    \qquad\quad u_l &= 1,
\\
    x_r &= 2, &
    \RHO_r &= 0.25, &
    u_r &= 0,
\end{alignat}
which gives an intermediate density $\RHO_m = 1.16641$, and shock
speeds $s_l = 0.36360$ and $s_r = 0.63640$. The exact velocity profile
can be computed as well, but we refer the reader to the literature.
Using a final time of $T=1.6$, we obtain
$$
\begin{tabular}{|c|c||c|c||c|c||c|c||c|c|}
    \hline
    \multicolumn{10}{|c|}{VPS1a}
\\
    \hline
    \text{$N$}
    & \text{$\tau$}
    & \text{$\WAS$}
    & \text{Rate}
    & \text{$E_\WAS$}
    & \text{Rate}
    & \text{$\L^1$}
    & \text{Rate}
    & \text{$E_\text{tot}$}
    & \text{Rate}
\\
    \hline
    \hline
    100   & .1    &  7.85e-3  &        &  3.39e-2  &        &  4.60e-2  &        &  1.79e-3  &
\\
    250   & .04   &  4.14e-3  &  0.70  &  2.14e-2  &  0.50  &  2.18e-2  &  0.82  &  7.62e-4  &  0.93
\\
    1000  & .01   &  1.55e-3  &  0.71  &  1.03e-2  &  0.53  &  6.10e-3  &  0.92  &  2.10e-4  &  0.93
\\
    2500  & .004  &  7.81e-4  &  0.75  &  5.99e-3  &  0.60  &  2.63e-3  &  0.92  &  8.87e-5  &  0.94
\\
    10000 & .0001 &  2.66e-4  &  0.78  &  2.07e-3  &  0.77  &  6.52e-4  &  1.01  &  2.33e-5  &  0.96
\\
    \hline
\end{tabular}
$$
$$
\begin{tabular}{|c|c||c|c||c|c||c|c||c|c|}
    \hline
    \multicolumn{10}{|c|}{VPS2}
\\
    \hline
    \text{$N$}
    & \text{$\tau$}
    & \text{$\WAS$}
    & \text{Rate}
    & \text{$E_\WAS$}
    & \text{Rate}
    & \text{$\L^1$}
    & \text{Rate}
    & \text{$E_\text{tot}$}
    & \text{Rate}
\\
    \hline
    \hline
    100    & .1    &  3.50e-3  &        &  2.28e-2  &        &  3.05e-2  &        &  8.05e-4  &
\\
    250    & .04   &  1.27e-3  &  1.11  &  1.60e-2  &  0.38  &  1.36e-2  &  0.88  &  3.77e-4  &  0.83
\\
    1000   & .01   &  2.89e-4  &  1.07  &  7.76e-3  &  0.52  &  4.48e-3  &  0.80  &  1.06e-4  &  0.92
\\
    2500   & .004  &  1.13e-4  &  1.02  &  3.67e-3  &  0.82  &  1.97e-3  &  0.90  &  3.38e-5  &  1.24
\\
    10000  & .0001 &  2.84e-5  &  1.00  &  1.29e-3  &  0.75  &  3.83e-4  &  1.18  &  7.94e-6  &  1.04
\\
    \hline
\end{tabular}
$$
Here $\WAS$ stands for the Wasserstein distance between the numerical
and the exact densities, while $\L^1$ indicates their
$\L^1$-difference. The functional $E_\WAS$ was defined in
\eqref{E:EW:DEF}, and we write $E_\mathrm{tot} :=
|E^n_\mathrm{tot}-E^\mathrm{exact} _\mathrm{tot}|$ for the difference
in the total energies. The $\L^\infty$-error of the densities does not
approach zero because of the shocks. We see that in spite of the high
frequency oscillations, the VPS2 scheme is more accurate than the VPS1a
method, with the density converging at first order and the velocity
converging somewhat slower. For both schemes, the total energy
$E_\text{tot}^n$ at the final time converges to
$E^\text{exact}_\text{tot}$ from above as $N\rightarrow\infty$ and
$\tau\rightarrow0$, thus the exact solution dissipates energy
(slightly) faster than any of the numerical solutions.

\begin{figure}
\includegraphics[scale=0.75]{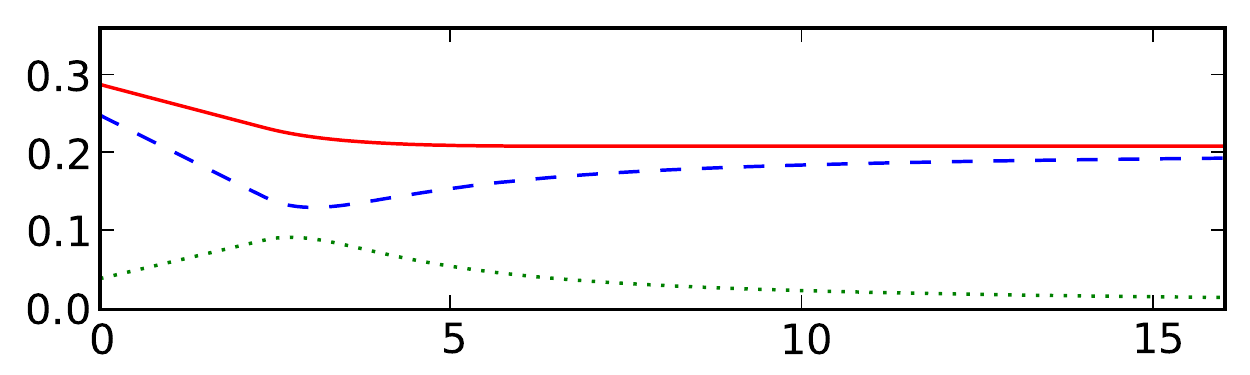}
\caption{Energy (solid: total, dashed: kinetic, dotted: internal).}
\label{F:SSENERGY}
\end{figure}

In Figure~\ref{F:SSENERGY}, we plot the energy of the numerical
solution (beyond the point that we are able to compute the exact
solution).
As long as the solution is discontinuous, the total energy decreases.
Moreover, the total energy is a linear function at the beginning of the
evolution (up to time $t=2$, say,) since the strength of the two shocks
remains constant. During this time, the internal energy increases along
with the width of the intermediate state $(\RHO_m,u_m)$, which carries
more energy because of its high density. As the shocks interact with
the rarefaction waves at the boundary of the support, their strengths
decrease and finally vanish (at about $t=6.7$), after which the total
energy remains essentially constant. This is in agreement with the
theory since for smooth solutions of the isentropic Euler equations,
the total energy is a conserved quantity.  Our schemes capture this
behavior remarkably well.


\subsubsection{Shock/Rarefaction}

\begin{figure}
\includegraphics[scale=0.75]{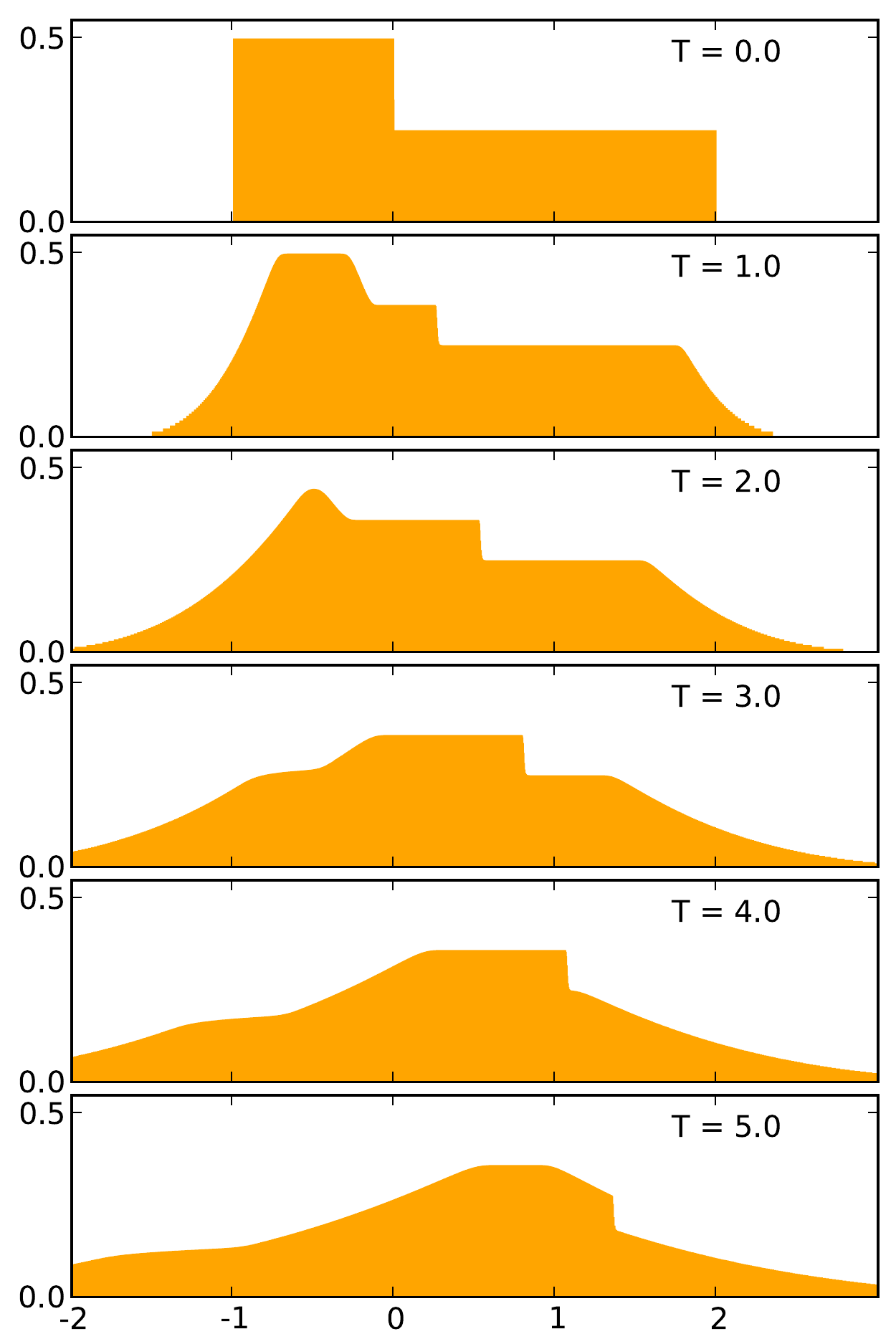}
\caption{Isentropic Euler equations with $\gamma=5/3$.}
\label{F:SRTIME}
\end{figure}

We also considered the case $\gamma=5/3$ with initial data
\begin{equation}
    (\bar{\RHO},\bar{u}) = \begin{cases}
        (0.5, 0) & \text{if $x\in(-1,0)$,}
\\
        (0.25, 0) & \text{if $x\in(0,2)$,}
\\
        (0, 0) & \text{otherwise.}
    \end{cases}
\label{E:SRDATA}
\end{equation}
The exact solution of the Riemann problem at $x=0$ is given by a
pattern involving a shock followed a rarefaction wave.
Figure~\ref{F:SRTIME}
shows the approximation at different times. The computation was done
using VPS1 with $m=0.001$ and $\tau=0.01$.

\begin{figure}
\includegraphics[scale=0.75]{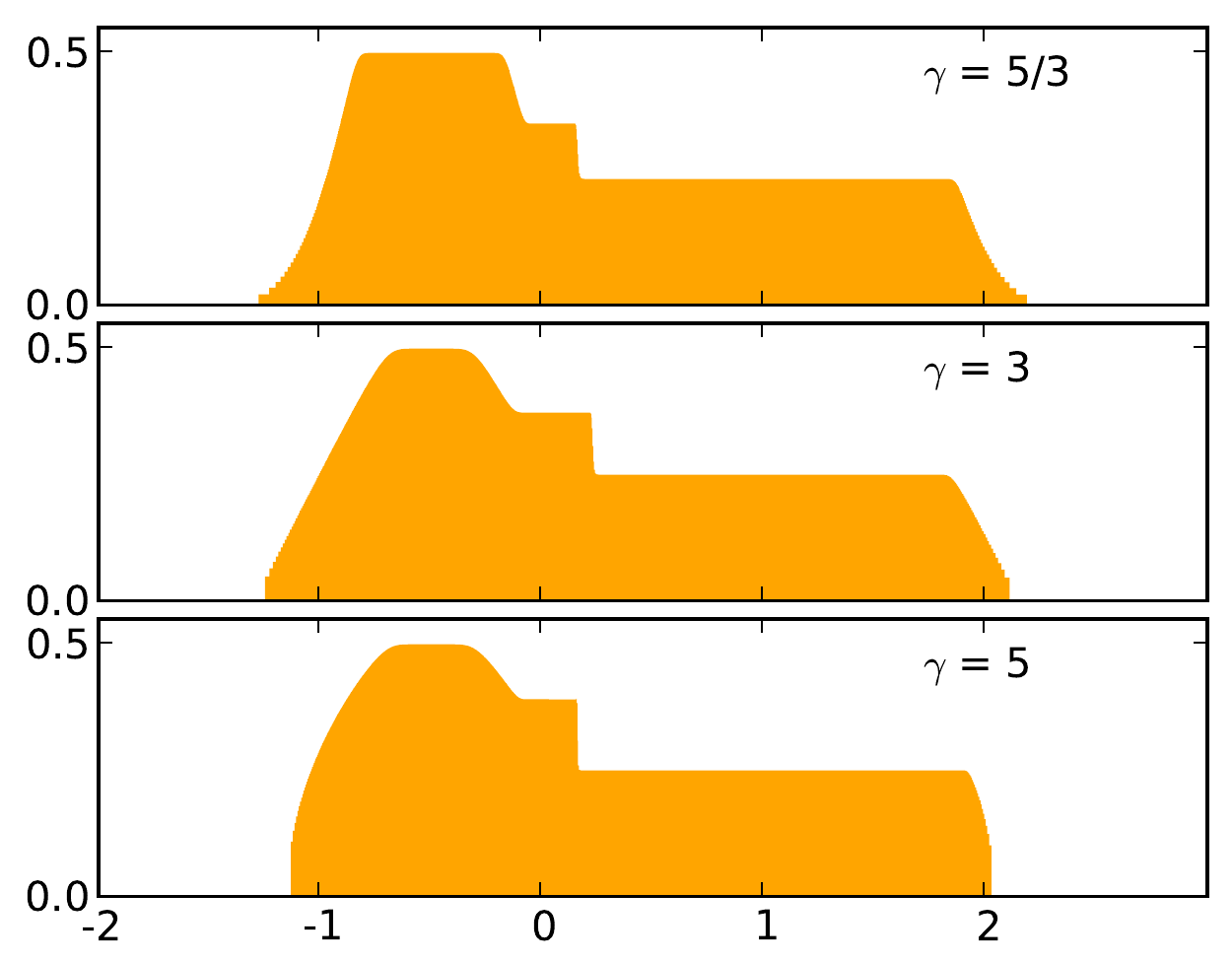}
\caption{Isentropic Euler equations at time $T=0.6$.}
\label{F:SRERROR}
\end{figure}

Figure~\ref{F:SRERROR}
shows the approximate solution of \eqref{E:IE} for the same parameters
and initial data \eqref{E:SRDATA}, at time $T=0.6$ for different
values of $\gamma$.  For sufficiently small $t>0$, the exact density
is given by the following formula:
\begin{equation}
    \RHO(t,x) = \begin{cases}
        \bigg( \DST\frac{(-u_l+\RHO_l^\theta) + (x-x_l)/t}{\theta+1}
                \bigg)^{1/\theta}
            & \text{if $x\in (x_1, x_2)$,}
\\
        \RHO_l
            & \text{if $x\in (x_2, x_3)$,}
\\
        \bigg( \DST\frac{(u_l+\RHO_l^\theta) - x/t}{\theta+1}
                \bigg)^{1/\theta}
            & \text{if $x\in (x_3, x_4)$,}
\\
        \RHO_m
            & \text{if $x\in (x_4, x_5)$,}
\\
        \RHO_r
            & \text{if $x\in (x_5, x_6)$,}
\\
        \bigg( \DST\frac{(u_r+\RHO_r^\theta) - (x-x_r)/t}{\theta+1}
                \bigg)^{1/\theta}
            & \text{if $x\in (x_6, x_7)$,}
\\
        0 & \text{otherwise,}
    \end{cases}
\label{E:PROFILE}
\end{equation}
with constants
\begin{align*}
    x_1 &= x_l + t (u_l-\RHO_l^\theta),
&
    x_3 &= t (u_l-\theta\RHO_l^\theta),
&
    x_6 &= x_r + t (u_r-\theta\RHO_r^\theta),
\\
    x_2 &= x_l + t (u_l+\theta\RHO_l^\theta),
&
    x_4 &= t (u_m-\theta\RHO_m^\theta),
&
    x_7 &= x_r + t (u_r+\RHO_r^\theta),
\\
&
&
    x_5 &= t s_r;
\end{align*}
see \cite{ChenWang}. For the case under consideration, we have
\begin{alignat*}{3}
    x_l &= -1, &
    \qquad\quad \RHO_l &= 0.5, &
    \qquad\quad u_l &= 0,
\\
    x_r &= 2, &
    \RHO_r &= 0.25, &
    u_r &= 0,
\end{alignat*}
which implies an intermediate density $\RHO_m = 0.3601$ and velocity
$u_m = 0.0823$, and a shock speed $s_r = 0.3636$. The exact velocity
profile can be computed as well, but we refer the reader to the
literature. The following table shows the $\L^1$-error:
$$
\begin{tabular}{|c||c|}
    \hline
    \text{Problem}
        & \text{$\L^1$-error}
\\
    \hline
    \hline
    $\gamma=5/3$
        & 4.46e-3 
\\
    \hline
    $\gamma=3$
        & 5.57e-3 
\\
    \hline
    $\gamma=5$
        & 3.87e-3 
\\
    \hline
\end{tabular}
$$
Notice that the rarefaction waves in \eqref{E:PROFILE} that connect the
profile to the vacuum, are convex for $\gamma<3$, linear if $\gamma=3$,
and concave if $\gamma>3$. To get a convergence rate, we computed the
approximate solution at time $t=0.6$ with initial data \eqref{E:SRDATA}
and $\gamma=5/3$ for different values of particle mass/timestep. We
have
$$
\begin{tabular}{|c|c||c|c|}
    \hline
    \text{$m$}
        & \text{$\tau$}
        & \text{$\L^1$-Error}
        & \text{Rate}
\\
    \hline
    \hline
    0.01
        & 0.1
        & 2.146e-2 
        &
\\
    \hline
    0.005
        & 0.05
        & 1.397e-2 
        & 0.620
\\
    \hline
    0.001
        & 0.01
        & 4.464e-3 
        & 0.709
\\
    \hline
    0.0005
        & 0.005
        & 2.615e-3 
        & 0.772
\\
    \hline
    0.0001
        & 0.001
        & 7.318e-4 
        & 0.791
\\
    \hline
\end{tabular}
$$
As with VPS1a in the previous section, the rate is somewhat less than
one.  We did not study the performance of the VPS1a or VPS2 schemes on
this initial data or implement the other measures of error in our VPS1
code.
%


\subsubsection{Rarefaction/Rarefaction}

\begin{figure}
\includegraphics[scale=0.75]{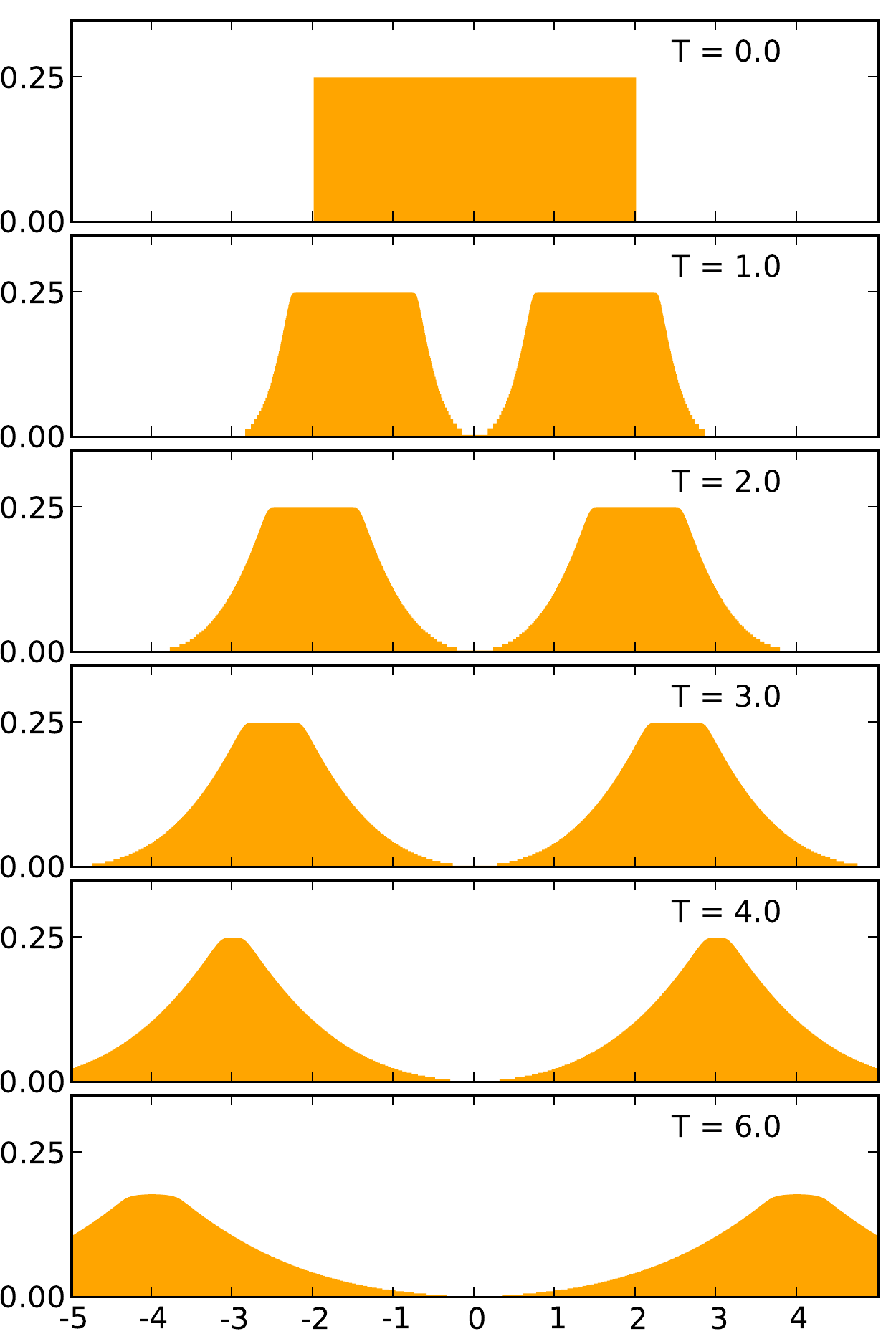}
\caption{Isentropic Euler equations with $\gamma=5/3$.}
\label{F:RRTIME}
\end{figure}

Finally, we use the VPS1 method to compute the approximate solution of
\eqref{E:IE} for $\gamma=5/3$ and initial data
$$
    (\bar{\RHO},\bar{u}) = \begin{cases}
        (0.25, -0.5) & \text{if $x\in(-2,0)$,}
\\
        (0.25, 0.5) & \text{if $x\in(0,2)$,}
\\
        (0, 0) & \text{otherwise.}
    \end{cases}
$$
Figure~\ref{F:RRTIME}
shows the result at different times. The fluid splits into two parts
that travel in opposite directions. Between the blocks, two rarefaction
waves form that are separated by vacuum. In our approximation, the
density is defined as particle mass divided by the specific volume,
which is the distance between neighboring particles. We therefore do
not get a perfect vacuum. As the two blocks move further apart,
however, the density between them does approach zero. By construction,
the density can never become negative.


\subsection*{Acknowledgements}

The research of M.~Westdickenberg was supported by NSF grant DMS
0701046.
The research of J.~Wilkening was supported in part by the Director,
Office of Science, Computational and Technology Research, U.S.
Department of Energy under Contract No. DE-AC02-05CH11231.


\appendix

\section{Porous Medium Equation}
\label{S:EU:LA}

We derive the Euler-Lagrange equations for the minimization problem
\begin{equation}
    \RHO^{n+1} := \ARGMIN\Bigg\{
        \frac{1}{2\tau} \WAS(\RHO^n,\RHO)^2 + \INT[\RHO]
            \colon \RHO\in\SP(\R^d) \Bigg\}
\label{E:DISGF2}
\end{equation}
to show that it indeed gives an approximation for the porous medium equation. More precisely, we will show that \eqref{E:DISGF2} can be interpreted as a backward Euler method applied to Darcy's law $\BU=-\nabla U'(\RHO)$ in a Lagrangian formulation of the problem. We recall first that if $\RHO^n$ and $\RHO$ are Lebesgue measurable functions, then there exists a Borel map $\BR:\R^d \longrightarrow \R^d$ defined $\RHO^n$-a.e., with the property that
\begin{equation}
    \WAS(\RHO^n,\RHO)^2 = \int_{\R^d} |\BR(x)-x|^2 \RHO^n(x) \,dx
    \quad\text{and}\quad
    \BR\#(\RHO^n\LEB^d) = \RHO\LEB^d.
\label{E:TE}
\end{equation}
The map $\BR$ is called an optimal transport map. It is the gradient
of a lower semicontinuous, convex function and invertible
$\RHO^n$-a.e. The second identity in \eqref{E:TE} now implies that for all test functions $\phi\in\CB(\R^d)$ we have
\begin{align*}
    \int_{\R^d} \varphi\big(\BR(x)\big) \RHO^n(x) \,dx
        &= \int_{\R^d} \varphi(z) \RHO(z) \,dz
            = \int_{\R^d} \varphi\big(\BR(x)\big) \RHO\big(\BR(x)\big)
                \det D\BR(x) \,dx.
\end{align*}
Here we used the change of variables formula, which can be justified because
$\BR$ is a monotone function. Since $\BR$ is invertible $\RHO^n$-a.e.\
we conclude that
$$
    \RHO^n(x) = \RHO\big(\BR(x)\big) \det D\BR(x)
    \quad\quad\text{for a.e.\ $x\in\R^d$,}
$$
which allows us to express the internal energy $\INT[\RHO]$ in the form
\begin{equation} \label{E:U:IN:TERMS:OF:R}
    \INT[\RHO] = \int_{\R^d} U\big(\RHO(z)\big) \,dz
        = \int_{\R^d} U\bigg( \frac{\RHO^n(x)}{\det D\BR(x)} \bigg)
            \det D\BR(x) \,dx.
\end{equation}
That is, we can consider the internal energy as a functional of $\BR$
instead of $\RHO$.

Let now $\RHO^{n+1}$ be the minimizer of \eqref{E:DISGF} and let
$\BR^{n+1}$ denote the optimal transport map that pushes $\RHO^n\LEB^d$
forward to the measure $\RHO^{n+1}\LEB^d$. For any smooth vector field
$\zeta\colon\R^d \longrightarrow\R^d$ and any $\EPS\GS 0$ we now define
the functions
$$
    \BR_\EPS := (\ID+\EPS\zeta)\circ\BR^{n+1}
    \quad\text{and}\quad
    \RHO_\EPS\LEB^d := \BR_\EPS \# (\RHO^n\LEB^d).
$$
Note that the map $\BR_\EPS$ is not an optimal transport map if
$\EPS>0$, but we have
$$
    \WAS(\RHO^n, \RHO_\EPS)^2
        \LS \int_{\R^d} |\BR_\EPS-\ID|^2 \RHO^n \,dx
    \quad\text{for all $\EPS\GS 0$.}
$$
This implies the estimate
\begin{align}
    \limsup_{\EPS\rightarrow 0}
        \frac{ \WAS(\RHO^n,\RHO_\EPS)^2
            -\WAS(\RHO^n,\RHO^{n+1})^2 }{\EPS}
        &\LS 2\int_{\R^d} (\zeta\circ\BR^{n+1})
            \cdot(\BR^{n+1}-\ID) \RHO^n \,dx
\label{E:ELFWAS}\\
    & = 2\int_{\R^d} \zeta\cdot\Big( (\BR^{n+1}-\ID )
        \circ (\BR^{n+1})^{-1} \Big) \RHO^{n+1} \,dz.
\nonumber
\end{align}
On the other hand, a straightforward computation using \eqref{E:PRESS} and
\eqref{E:U:IN:TERMS:OF:R} shows that
\begin{align}
\notag
    \lim_{\EPS\rightarrow 0} \frac{\INT[\RHO_\EPS]
            -\INT[\RHO^{n+1}]}{\EPS}
        &= \int_{\R^d} \Bigg\{
            -U'\bigg( \frac{\RHO^n}{\det D\BR^{n+1}} \bigg)
                \frac{\RHO^n}{\det D\BR^{n+1}}
            +U\bigg( \frac{\RHO^n}{\det D\BR^{n+1}} \bigg) \Bigg\}
\\
\notag
        & \quad\qquad \times
            \,\mathrm{tr}\Big( (D\zeta)\circ\BR^{n+1} \Big)
            \,\det D\BR^{n+1} \,dx
\\
\label{E:VAR:DER:U}
        &= -\int_{\R^d} P(\RHO^{n+1}) \,\nabla\cdot\zeta \,dz.
\end{align}
We used the change of variables formula again. Now \eqref{E:DISGF2} implies that
$$
    0 \LS \frac{1}{\tau}\int_{\R^d} \zeta\cdot\Big( (\BR^{n+1}-\ID )
            \circ (\BR^{n+1})^{-1} \Big) \RHO^{n+1} \,dz
        -\int_{\R^d} P(\RHO^{n+1}) \,\nabla\cdot\zeta \,dz
$$
for all test functions $\zeta$. Note that by changing the sign of $\zeta$ we obtain the converse inequality, so we actually have equality here. Since $\zeta$ was arbitrary we obtain
$$
    \bigg( \frac{\BR^{n+1}-\ID}{\tau} \circ (\BR^{n+1})^{-1} \bigg)
        \RHO^{n+1} + \nabla P(\RHO^{n+1}) = 0
$$
in the sense of distributions. If we define $\BU^{n+1} :=
\tau^{-1}(\BR^{n+1}-\ID) \circ (\BR^{n+1})^{-1}$, then we obtain
Darcy's law $\BU^{n+1} = -\nabla U'(\RHO^{n+1})$ a.e.\ in $\{\RHO^{n+1}
>0\}$. Since the new density $\RHO^{n+1}$ can be computed according to
the formula $\RHO^n = (\RHO^{n+1}\circ\BR^{n+1}) \,\det D\BR^{n+1}$
a.e., we conclude that the minimization problem \eqref{E:DISGF} in fact
amounts to a backward Euler method for the transport map $\BR^{n+1}$
pushing $\RHO^n\LEB^d$ forward to $\RHO^{n+1}\LEB^d$.


\section{Isentropic Euler Equations}
\label{S:DERIVE:VTD}

The time discretization of Definition~\ref{D:VTD} is motivated by
Dafermos' entropy rate criterion \cite{Dafermos}: dissipate the total
energy as fast as possible. This objective is achieved through two
different mechanisms. First, the kinetic energy is minimized by
replacing the given velocity by the corresponding optimal transport
velocity. Second, the internal energy is minimized by solving an
optimization problem similar to the one in \eqref{E:DISGF} for the
porous medium equation, but with a different penalty term. The
intuition for the latter is the following: We think of the fluid as a
collection of particles, each one determined by a position and a
velocity. The particles prefer to stay on their free flow trajectories
(minimizing the acceleration), but the pressure from the surrounding
particles forces them to deviate from their characteristic paths.

To measure the acceleration we introduce a functional similar to the
Wasserstein distance. Consider first a single particle that initially
is located at position $x_1\in\R^d$ with initial velocity
$\xi_1\in\R^d$. Specify some timestep $\tau>0$. During the time
interval of length $\tau$ we allow the particle to move to a new
position $x_2\in\R^d$ and to change its velocity to $\xi_2\in\R^d$.
What is the minimal acceleration required to achieve this change? We
are looking for a curve $c\colon[0,\tau]\longrightarrow\R^d$ such that
$$
    (c,\dot{c})(0) = (x_1,\xi_1)
    \quad\text{and}\quad
    (c,\dot{c})(\tau) = (x_2,\xi_2),
$$
that minimizes the $\L^2([0,\tau])$-norm of the second derivative
$\ddot{c}$ along the curve. This minimization problem has a unique
solution; the minimizer is given by a cubic polynomial; and the minimal
(averaged) acceleration is given by the formula
$$
    \frac{1}{\tau} \int_0^\tau |\ddot{c}(s)|^2 \,ds
        = 12 \bigg| \frac{1}{\tau} \bigg( \frac{x_2-x_1}{\tau}
            -\frac{\xi_2+\xi_1}{2} \bigg) \bigg|^2
        + \bigg| \frac{\xi_2-\xi_1}{\tau} \bigg|^2.
$$
%

Consider now a fluid characterized by a density $\RHO$ and a velocity
field $\BU$. Assume that we are given a transport map $\BR\colon\R^d
\longrightarrow\R^d$ that is invertible $\RHO$-a.e.\ and determines
where each particle travels to during a time interval of length
$\tau>0$. Which velocity field $\hat{\BU}$ minimizes the total
acceleration cost, defined as the integral
\begin{equation}
    \mathscr{A}_{\tau,\BR}[\hat{\BU}]^2
        := \int_{\R^d} \Bigg( 3 \bigg| \frac{\BR-\ID}{\tau}
                -\frac{\hat{\BU}\circ\BR+\BU}{2} \bigg|^2
            + \frac{1}{4} |\hat{\BU}\circ\BR-\BU|^2 \Bigg) \RHO \,dx?
\label{E:ATR}
\end{equation}
Note that since $\BR$ is assumed to be invertible, there is only a
single particle located at position $\BR(x)$ for $\RHO$-a.e.\
$x\in\R^d$, and so the velocity of the transported particles can be
described by an Eulerian velocity field $\hat{\BU}$. Let $\BU^+$ denote
the minimizer of the functional $\mathscr{A}_{\tau,\BR}$ (we refer the
reader to \cite{GangboWestdickenberg} for further details). For any
smooth vector field $\zeta\colon\R^d \longrightarrow\R^d$ and $\EPS\GS
0$ define $\BU_\EPS := \BU^+ + \EPS\zeta$. Then we compute
$$
    \lim_{\EPS\rightarrow 0}
        \frac{\mathscr{A}_{\tau,\BR}[\BU_\EPS]^2
            -\mathscr{A}_{\tau,\BR}[\BU^+]^2}{\EPS}
        = \int_{\R^d} \bigg( -3\frac{\BR-\ID}{\tau}
            + 2\BU^+\circ\BR + \BU \bigg) \cdot (\zeta\circ\BR)
                \RHO \,dx,
$$
which vanishes because $\BU^+$ is a minimizer. Since $\zeta$ was
arbitrary, we obtain that
\begin{equation}
    \BU^+\circ\BR = \BU + \frac{3}{2\tau} \Big( \BR-(\ID+\tau\BU) \Big)
\label{E:VELOI}
\end{equation}
$\RHO$-a.e. Using this identity in \eqref{E:ATR} we have
$$
    \mathscr{A}_{\tau,\BR}[\BU^+]^2
        = \frac{3}{4\tau^2} \int_{\R^d} |\BR-(\ID+\tau\BU)|^2 \RHO \,dx
        =: \mathscr{W}_\tau[\BR]^2.
$$
Assume now that the push-forward of $\RHO\LEB^d$ under the map
$\ID+\tau\BU$ is absolutely continuous with respect to the Lebesgue
measure and define $\hat{\RHO}\LEB^d := (\ID+\tau\BU)\#(\RHO\LEB^d)$.
If we fix a density $\RHO^+\in\SP(\R^d)$, then there are many transport
maps $\BR\colon\R^d\longrightarrow\R^d$ that push $\RHO\LEB^d$ forward
to the measure $\RHO^+\LEB^d$. Minimizing in $\BR$, we obtain
$$
    \inf\Big\{ \mathscr{W}_\tau[\BR]^2 \colon
        \BR\#(\RHO\LEB^d) = \RHO^+\LEB^d \Big\}
        = \frac{3}{4\tau^2} \WAS(\hat{\RHO},\RHO^+)^2,
$$
so the minimal total acceleration cost is proportional to the
Wasserstein distance between $\RHO^+$ and the density $\hat{\RHO}$
obtained by freely transporting the initial density $\RHO$ along the
initial velocity field $\BU$. We refer again to
\cite{GangboWestdickenberg} for further details.

Note that Step~2 of the time discretization in Definition~\ref{D:VTD}
amounts to minimizing the internal energy $\INT[\RHO]$ while penalizing
the deviation of particles from their characteristic trajectories
(which requires accelerating the particles). Computing the
Euler-Lagrange equations for the minimization problem \eqref{E:DISGF},
as we did above for the time step for the porous medium equation, we
find that the map
$$
    \BT^{n+1} := \ID + \frac{2\tau^2}{3} \nabla U'(\RHO^{n+1})
$$
is an optimal transport map between $\RHO^{n+1}\LEB^d$ and
$\hat{\RHO}^n\LEB^d$, thus invertible $\RHO^{n+1}$-a.e.

The update of the velocity \eqref{E:UN:UP} follows from formula
\eqref{E:VELOI} if we set
$$
    \BU := \hat{\BU}^n,
    \quad
    \BU^+ := \BU^{n+1},
    \quad
    \BR := \bigg( \ID + \frac{2\tau^2}{3} \nabla U'(\RHO^{n+1})
        \bigg)^{-1} \circ (\ID+\tau\BU^n).
$$
Indeed we have that
\begin{align*}
    & \frac{3}{2\tau} \Bigg( \bigg( \ID+\frac{2\tau^2}{3} \nabla
        U'(\RHO^{n+1}) \bigg)^{-1} - \ID \Bigg) \circ \bigg(
            \ID+\frac{2\tau^2}{3} \nabla U'(\RHO^{n+1}) \bigg)
\\
    & \quad
        = \frac{3}{2\tau} \Bigg( \ID - \bigg( \ID+\frac{2\tau^2}{3}
            \nabla U'(\RHO^{n+1}) \bigg) \Bigg)
        = -\tau\nabla U'(\RHO^{n+1})
\end{align*}
$\RHO^{n+1}$-a.e. The only caveat is that \eqref{E:VELOI} was derived
under the assumption that the map $\BR$ is essentially invertible. But
this condition is satisfied because $\hat{\BU}^n$ was chosen in such a
way that the map $\ID+\tau\hat{\BU}^n$ is an optimal transport map
pushing $\RHO^n\LEB^d$ forward to $\hat{\RHO}^n \LEB^d :=
(\ID+\tau\BU^n)\# (\RHO^n\LEB^d)$; see Step~1 in
Definition~\ref{D:VTD}. As mentioned there, this step decreases the
kinetic energy in an optimal way.



\begin{bibdiv}
\begin{biblist}

\bib{AmbrosioGigliSavare}{book}{
    AUTHOR = {Ambrosio, L.},
    AUTHOR = {Gigli, N.},
    AUTHOR = {Savar{\'e}, G.},
     TITLE = {Gradient flows in metric spaces and in the space of
              probability measures},
    SERIES = {Lectures in Mathematics ETH Z\"urich},
 PUBLISHER = {Birkh\"auser Verlag},
   ADDRESS = {Basel},
      YEAR = {2005},
}

\bib{ArnoldKhesin}{book}{
    AUTHOR = {Arnold, V. I.},
    AUTHOR = {Khesin, B. A.},
     TITLE = {Topological methods in hydrodynamics},
    SERIES = {Applied Mathematical Sciences},
    VOLUME = {125},
 PUBLISHER = {Springer-Verlag},
   ADDRESS = {New York},
      YEAR = {1998},
     PAGES = {xvi+374},
      ISBN = {0-387-94947-X},
}

\bib{Caffarelli}{incollection}{
   AUTHOR = {L. A. Caffarelli},
    TITLE = {Allocation maps with general cost functions},
BOOKTITLE = {Partial differential equations and applications},
   EDITOR = {P. Marcellini and G. G. Talenti and E. Vesintini},
   SERIES = {Lecture Notes in Pure and Applied Mathematics},
   VOLUME = {177},
PUBLISHER = {Marcel Dekker, Inc},
  ADDRESS = {New York},
    PAGES = {29--35},
     YEAR = {1996},
}

\bib{ChenWang}{incollection}{
    AUTHOR = {Chen, G.-Q.},
    AUTHOR = {Wang, D.},
     TITLE = {The {C}auchy problem for the {E}uler equations for
              compressible fluids},
 BOOKTITLE = {Handbook of mathematical fluid dynamics, Vol. I},
     PAGES = {421--543},
 PUBLISHER = {North-Holland},
   ADDRESS = {Amsterdam},
      YEAR = {2002},
}

\bib{cvxopt}{misc}{
      NOTE = {\url{http://abel.ee.ucla.edu/cvxopt}},
}

\bib{Dafermos}{article}{
    AUTHOR = {Dafermos, C. M.},
     TITLE = {The entropy rate admissibility criterion for solutions of
              hyperbolic conservation laws},
   JOURNAL = {J. Differential Equations},
    VOLUME = {14},
      YEAR = {1973},
     PAGES = {202--212},
}

\bib{GangboMcCann}{article}{
    AUTHOR = {W. Gangbo and R. J. McCann},
     TITLE = {The geometry of optimal transportation},
   JOURNAL = {Acta Math.},
    VOLUME = {177},
    NUMBER = {2},
     PAGES = {113--161},
      YEAR = {1996},
}

\bib{GangboWestdickenberg}{article}{
    AUTHOR = {Gangbo, W.},
    AUTHOR = {Westdickenberg, M.},
     TITLE = {Optimal Transport for the system of Isentropic Euler
              equations},
      NOTE = {FIXME},
      YEAR = {2008},
}

\bib{HairerNorsettWanner}{book}{
    AUTHOR = {Hairer, E.},
    AUTHOR = {Norsett, S. P.},
    AUTHOR = {Wanner, G.},
     TITLE = {Solving Ordinary Differential Equations {I}:
              Nonstiff Problems},
 PUBLISHER = {Springer},
   ADDRESS = {Berlin},
   EDITION = {2nd},
      YEAR = {2000}
}

\bib{HolmMarsdenRatiu}{article}{
    AUTHOR = {Holm, D. D.},
    AUTHOR = {Marsden, J. E.},
    AUTHOR = {Ratiu, T. S.},
     TITLE = {The {E}uler-{P}oincar\'e equations and semidirect products
              with applications to continuum theories},
   JOURNAL = {Adv. Math.},
    VOLUME = {137},
      YEAR = {1998},
    NUMBER = {1},
     PAGES = {1--81},
}

\bib{KinderlehrerWalkington}{article}{
    AUTHOR = {Kinderlehrer, D.},
    AUTHOR = {Walkington, N. J.},
     TITLE = {Approximation of parabolic equations using the {W}asserstein
              metric},
   JOURNAL = {M2AN Math. Model. Numer. Anal.},
    VOLUME = {33},
      YEAR = {1999},
    NUMBER = {4},
     PAGES = {837--852},
}

\bib{Knitro}{misc}{
      NOTE = {\url{http://www.ziena.com/knitro.htm}},
}

\bib{MarsdenWest}{article}{
    AUTHOR = {Marsden, J. E.},
    AUTHOR = {West, M.},
     TITLE = {Discrete mechanics and variational integrators},
   JOURNAL = {Acta Numer.},
    VOLUME = {10},
      YEAR = {2001},
     PAGES = {357--514},
}
	
\bib{NocedalWright}{book}{
    AUTHOR = {Nocedal, J.},
    AUTHOR = {Wright, S. J.},
     TITLE = {Numerical Optimization},
 PUBLISHER = {Springer},
   ADDRESS = {New York},
      YEAR = {1999}
}

\bib{Otto}{article}{
    AUTHOR = {Otto, F.},
     TITLE = {The geometry of dissipative evolution equations: the
              porous medium equation},
   JOURNAL = {Comm. Partial Differential Equations},
    VOLUME = {26},
      YEAR = {2001},
    NUMBER = {1-2},
     PAGES = {101--174},
}

\bib{Vazquez}{article}{
    AUTHOR = {V{\'a}zquez, J. L.},
     TITLE = {Perspectives in nonlinear diffusion: between
              analysis, physics and geometry},
 BOOKTITLE = {International Congress of Mathematicians. Vol. I},
     PAGES = {609--634},
 PUBLISHER = {Eur. Math. Soc., Z\"urich},
      YEAR = {2007},
}

\bib{Villani}{book}{
    AUTHOR = {Villani, C.},
     TITLE = {Topics in optimal transportation},
    SERIES = {Graduate Studies in Mathematics},
    VOLUME = {58},
 PUBLISHER = {American Mathematical Society},
   ADDRESS = {Providence, R.I.},
      YEAR = {2003},
}



\end{biblist}
\end{bibdiv}
\end{document}